\documentclass{article}

\usepackage{delarray,verbatim,enumerate,a4wide}
\usepackage{amsmath,amsthm,amstext,amsbsy,amssymb,amsfonts,amscd}
\usepackage[all]{xy}
\usepackage[frenchb]{babel}
\usepackage[T1]{fontenc}
\usepackage[utf8]{inputenc}

\usepackage{scalerel}
\usepackage[small]{titlesec}

\usepackage{color}
\definecolor{vert}{rgb}{0.1,0.4,0.2}
\usepackage[colorlinks=true,linkcolor=blue,citecolor=vert,linktocpage]{hyperref}

\usepackage{calligra}
\DeclareFontShape{T1}{calligra}{m}{n}{<->s*[0.95]callig15}{}
\DeclareMathAlphabet{\mathscr}{T1}{calligra}{m}{n}

\newtheorem{Th}{Théorème}[]
\newtheorem{Lem}[Th]{Lemme}
\newtheorem{Prop}[Th]{Proposition}
\newtheorem{Cor}[Th]{Corollaire}

\newtheorem{Sco}[Th]{Scolie}
\newtheorem{Def} [Th]{Définition}

\newtheorem*{ThA}{Théorème A}
\newtheorem*{ThB}{Théorème B}
\newtheorem*{ThC}{Théorème C}

\def\Preuve{\noindent {\it Preuve.~}}

\def\Remarque{\smallskip\noindent {\it Remarque.~}}

\def\Exemple{\smallskip\noindent {\it Exemple.~}}

\font\teneufm=eufm10
\font\seveneufm=eufm7
\font\fiveeufm=eufm5
\newfam\gothfam
\textfont\gothfam=\teneufm \scriptfont\gothfam=\seveneufm
\scriptscriptfont\gothfam=\fiveeufm

		\def\QQ{\mathbb Q}	
\def\NN{\mathbb N}	\def\ZZ{\mathbb Z}	\def\FF{\mathbb F}	
\def\F2{\mathbb{F}_2}	\def\Z2{\mathbb{Z}_2}		
\def\Zl{\mathbb{Z}_\ell} 	\def\Ql{\mathbb{Q}_\ell}		

 		\def\P{\mathcal  P}		\def\U{\mathcal  U}		\def\F{\mathcal  F}
\def\J{\mathcal  J}  		\def\C{\mathcal  C}		\def\R{\mathcal  R}		\def\X{\mathcal  X}
\def\Dl{\mathcal  D\ell} 	\def\Pl{\mathcal  P\ell}  	\def\Cl{\mathcal  C\!\ell}	\def\V{\mathcal  V}
\def\E{\mathcal  E}		\def\T{\mathcal  T}		\def\Dl{\mathcal  D\ell}	\def\D{\mathcal  D}

	\def\p{{\mathfrak p}}		\def\a{{\mathfrak a}}	
		\def\l{{\mathfrak l}}

\def\wi{\widetilde}		
	
	\def\deg{\operatorname{deg}}
\def\Gal{\operatorname{Gal}}		
\def\Ker{\operatorname{Ker}}		

\newcommand\scale[2]{\vstretch{#1}{\hstretch{#1}{#2}}}

\makeatletter
\newcommand*\wt[2][0.2ex]{%
        \begingroup
        \mathchoice{\wt@helper{#1}{#2}{\displaystyle}{\textfont}}
                   {\wt@helper{#1}{#2}{\textstyle}{\textfont}}
                   {\wt@helper{#1}{#2}{\scriptstyle}{\scriptfont}}
                   {\wt@helper{#1}{#2}{\scriptscriptstyle}{\scriptscriptfont}}%
        \endgroup
        #2%
}
\newcommand*\wt@helper[4]{%
        \def\currentfont{\the#41}%
        \def\currentskewchar{\char\the\skewchar\currentfont}%
        \setbox\tw@\hbox{\currentfont$#2$\currentskewchar}%
        \dimen@ii\wd\tw@
        \setbox\tw@\hbox{\currentfont$#2${}\currentskewchar}%
        \advance\dimen@ii-\wd\tw@
        \rlap{\raisebox{-#1}{$\m@th#3\kern\dimen@ii\widetilde{\phantom{#2}}$}}%
}
\makeatother

\def\wE{\,\wt[0.1ex]{\!\mathcal E}}		\def\wU{\wt[0.2ex]{\mathcal U}}
\def\wJ{\,\wt[0.2ex]{\!\mathcal J}}	\def\wCl{\wt[0.1ex]{\mathcal C\!\ell}} \def\wDl{\wt[0.2ex]{\mathcal D\!\ell}}
	\def\wdiv{\wt[0.3ex]{div}}			\def\wnu{\wt[0.1ex]{\nu}}
\def\wCap{\wt[0.3ex]{Cap}}
\begin{document}

\title{\Large\bf Note sur la conjecture de Greenberg\footnote{J. Ramanujan Math. Soc. (à paraître)   }}

\author{ Jean-François {\sc Jaulent} }
\date{}
\maketitle
\bigskip

{\small
\noindent{\bf Résumé.} Nous étudions la conjecture de Greenberg sur les $\ell$-invariants d'Iwasawa des corps totalement réels, en termes de classes logarithmiques. Sous la conjecture de Leopoldt, et donc en particulier pour les corps abéliens, nous montrons inconditionnellement que la conjecture de Greenberg postule tout simplement la capitulation des classes logarithmiques dans la $\Zl$-tour cyclotomique. Dans le cas complètement décomposé, nous prouvons de même que le critère suffisant introduit par Gras exprime simplement la trivialité du $\ell$-groupe des classes logarithmiques du corps de base considéré. Enfin, dans le cas abélien, nous explicitons la structure des $\ell$-groupes de classes circulaires dans le contexte de la conjecture faible.}

\

{\small
\noindent{\bf Abstract.} We use logarithmic $\ell$-class groups to take a new view on Greenberg's conjecture about Iwasawa $\ell$-invariants of a totally real number field $K$. By the way we recall and complete some classical results. Under Leopoldt's conjecture, we unconditionally prove that Greenberg's conjecture holds if and only if the logarithmic classes of $K$ principalize in the cyclotomic $\Zl$-extension. As an illustration of our approach, in the special case where the prime $\ell$ splits completely in $K$, we prove that the sufficient condition introduced by Gras just asserts the triviality of the logarithmic class group of $K$.  Last, in the abelian case, we provide an explicit description of the circular class groups in connexion with the so-called weak conjecture.}
\bigskip

\tableofcontents
%%%%%%%%%%%%%%%%%%%%%%%%%%%%%%%%%%%%%%%%%%%%%%%%%%%%%%%%
%%%%%%%%%%%%%%%%%%%%%%%%%%%%%%%%%%%%%%%%%%%%%%%%%%%%%%%%
\section*{Introduction}
\addcontentsline{toc}{section}{Introduction}
\medskip
%%%%%%%%%%%%%%%%%%%%%%%%%%%%%%%%%%%%%%%%%%%%%%%%%%%%%%%%
%%%%%%%%%%%%%%%%%%%%%%%%%%%%%%%%%%%%%%%%%%%%%%%%%%%%%%%%

Le résultat emblématique de la Théorie d'Iwasawa affirme que les ordres respectifs $\ell^{h_n}$ des $\ell$-sous-groupe de Sylow $\,\Cl_{K_n}$ des groupes de classes d'idéaux attachés aux étages finis $K_n$ d'une $\Zl$-extension $K_\infty=\cup_{n\in\NN}K_n$ d'un corps de nombres $K$ sont donnés asymptotiquement par une formule de la forme:

\centerline{$h_n=\mu^{\phantom{*}}_K\ell^n + \lambda^{\phantom{*}}_K n + \nu^{\phantom{*}}_K$,}\smallskip

\noindent où $\lambda_K^{\phantom{*}}$ et $\mu_K^{\phantom{*}}$ sont des invariants structurels entiers attachés à la limite projective pour les morphismes normiques $\,\C_K=\varprojlim\,\Cl_{K_n}$ et $\nu_K^{\phantom{*}}$ un entier relatif convenable.\smallskip

Il est conjecturé, et cela a été démontré pour $K$ abélien par Ferrero et Washington, que l'invariant $\mu_K^{\phantom{*}}$ est toujours nul lorsque $K_\infty$ est la $\Zl$-extension cyclotomique de $K$.\smallskip

Lorsque, en outre, le corps de base $K$ est totalement réel, la conjecture de Leopoldt (qui est vérifiée elle aussi pour $K$ abélien en vertu du théorème d'indépendance de logarithmes de Baker-Brumer) affirme que $K$ ne possède d'autre $\Zl$-extension que la cyclotomique, qui provient, par composition avec $K$, de l'unique pro-$\ell$-extension abélienne réelle $\QQ_\infty$ de $\QQ$ qui est $\ell$-ramifiée (i.e. non ramifiée en dehors de $\ell$).\smallskip

Dans ce même cadre, i.e. pour $K$ totalement réel, la conjecture de Greenberg postule que les invariants structurels $\lambda_K^{\phantom{*}}$ et $\mu_K^{\phantom{*}}$ sont tous deux nuls,; autrement dit que les $\ell$-groupes de classes $\,\Cl_{K_n}$ sont stationnaires lorsqu'on monte la tour; ce qui revient à affirmer que leur limite projective $\,\C_K$ est un groupe fini.\smallskip

Tester numériquement cette conjecture s'avère rapidement problématique, même pour un corps totalement réel de petit degré $d$, car elle est par essence de nature asymptotique, ce qui nécessite, en principe, d'effectuer des calculs dans des étages de la $\Zl$-tour de plus en plus élevés donc de degrés $[K_n:\QQ]=d\ell^n$ de plus en plus hauts. D'où l'intérêt d'obtenir des critères de validité susceptibles de se lire, sinon directement dans $K$, du moins à des étages assez bas de la tour.
\medskip

Nous nous appuyons pour cela sur la notion de classe logarithmique qui, de par sa définition même, se trouve particulièrement adaptée à l'étude des phénomènes cyclotomiques. Il n'entre naturellement pas dans le cadre de ce travail de faire une présentation complète de l'arithmétique des classes et unités logarithmiques: nous nous appuyons sur \cite{J28} pour leur définition; sur \cite{J31} pour les éléments de la théorie $\ell$-adique du corps de classes qui sous-tendent leurs principales propriétés; sur \cite{BJ} et \cite{DJ+} enfin pour une description des algorithmes qui en permettent le calcul.\smallskip

Disons simplement ici que ces groupes s'obtiennent naturellement en envoyant le $\ell$-adifié $\R_K=\Zl\otimes_\ZZ K^\times$ du groupe multiplicatif du corps considéré dans la somme formelle $\,\Dl_K=\bigoplus_\p \Zl\,\p$ construite sur les places finies de $K$ par une famille de valuations $\wdiv=(\wnu_\p)_\p$ qui coïncident avec les valuations usuelles $(\nu_\p)_\p$ à l'exception des places au-dessus de $\ell$ pour lesquelles leur définition fait intervenir le logarithme $\ell$-adique de la norme locale. Dans la suite exacte obtenue:\smallskip

\centerline{$\begin{CD}1 @>>> \wE_K @>>> \Zl\otimes_\ZZ K^\times @>\wdiv>> \wi{\bigoplus}_\p\; \Zl\,\p @>>> \wCl_K @>>> 1,\end{CD}$}\smallskip

\noindent restreinte aux diviseurs de degré nul $\wDl_K=\wi{\bigoplus}_\p\; \Zl\,\p = \{\a=\sum_\p \nu_\p\,\p\;|\; \deg \a =\sum_\p \nu_\p\,\deg\p=0\}$,
les groupes de classes $\,\wCl_K$ et d'unités logarithmiques $\wE_K$ apparaissent respectivement comme noyau et conoyau du morphisme $\wdiv$
et donc comme analogues des $\ell$-adifiés des groupes de classes $\,\Cl_K$ et d'unités $\E_K=\Zl\otimes_\ZZ E_K$ au sens ordinaire.\smallskip

En dehors du cas abélien (où cela résulte encore du théorème d'indépendance de logarithmes de Baker-Brumer) et de quelques autres (cf. e.g. \cite{J38}), on ne sait pas si le $\ell$-groupe des classes logarithmiques est fini: de fait, comme expliqué dans \cite{J31,J55}, la finitude de  $\,\wCl_K$ est équivalente à la conjecture de Gross-Kuz'min pour le corps $K$ et le premier $\ell$, initialement énoncée dans \cite{Kuz}.

Néanmoins cet obstacle théorique ne s'oppose pas au calcul: dès lors qu'on connait le corps $K$ numériquement en ce sens qu'on est capable de déterminer ses invariants arithmétiques, il est toujours possible de calculer le groupe des classes logarithmiques et d'en tester la finitude comme la trivialité. Et, pour $K$ totalement réel donné, on peut même espérer que  $\,\wCl_K$  soit trivial pour presque tout $\ell$. \medskip

Ce point précisé, les résultats théoriques fondamentaux de ce travail sont les suivants:\medskip

-- En premier lieu nous montrons que le critère de capitulation donné par Greenberg dans le cas où le corps totalement réel $K$ ne possède qu'une seule place au-dessus de $\ell$ (cf. \cite{Grb}, Th. 1) se généralise très simplement sous la conjecture de Leopoldt en termes de classes logarithmiques quelle que soit la décomposition de $\ell$, comme observé indépendamment par Nguyen Quang Do (cf. \cite{Ng3}, Th. C) en liaison avec une version antérieure de ce travail:

\begin{ThA}[Théorème \ref{CGLog} \& Scolie] Si $K$ est un corps de nombres abélien totalement réel et $\ell$ un nombre premier arbitraire, la conjecture de Greenberg est vérifiée en $\ell$ dans $K$ si et seulement si le $\ell$-groupe des classes logarithmiques $\,\wCl_K$ capitule à un étage fini de la $\Zl$-tour cyclotomique $K_\infty/K$, dès lors que la conjecture de Leopoldt est satisfaite dans $K_\infty$ (par exemple pour $K$ abélien).
\end{ThA}
\smallskip

-- En second lieu nous montrons que le critère de validité étudié par Gras (cf. \cite{Gra5}, Th. 3.4) dans le cas complètement décomposé en liaison avec les heuristiques sur les régulateurs $\ell$-adiques (cf. \cite{Gra4}) s'interprète en termes de classes logarithmiques:

\begin{ThB}[Théorème \ref{CSGras} \& Scolie]
Si $K$ est un corps de nombres totalement réel qui vérifie la conjecture de Leopoldt pour un premier donné $\ell$, la conjecture de Greenberg est vérifiée en $\ell$ dans $K$ sous la condition suffisante que le $\ell$-groupe des classes logarithmiques $\,\wCl_K$ soit trivial. De plus, lorsque $\ell$ est complètement décomposé dans $K$, cette condition coïncide avec le critère de Gras.
\end{ThB}
\smallskip

-- Dans une dernière partie enfin, nous nous intéressons à la trivialité du sous-module fini $\,\F_K$ du $\Lambda$-module $\,\C_K$ en liaison avec ce qu'il est convenu d'appeler la {\em conjecture faible}, i.e. l'implication:\smallskip

\centerline{$\,\F_K=1 \Rightarrow \,\C_K=1$.}\smallskip

\noindent Nous montrons en particulier que l'hypothèse de trivialité $\,\F_K=1$, qui s'interprète classiquement en termes de capitulation, se lit directement dans le cas décomposé sur les sous-groupes sauvages, i.e. sur les sous-groupes des $\ell$-groupes de classes d'idéaux construits sur les places au-dessus de $\ell$ aux divers étages finis $K_n$ de la $\Zl$-tour cyclotomique $K_\infty/K$:

\begin{ThC}[Théorème \ref{CrSauvage}]
Soient $K$ un corps totalement réel et $\ell$ un premier complètement décomposé dans $K$. Sous la conjecture de Leopodt pour $\ell$ dans $K$, on a l'équivalence:\smallskip

\centerline{$\,\F^{\phantom{*}}_K=1 \quad\Leftrightarrow\quad \C^{[\ell]}_K=\varprojlim \,\Cl^{[\ell]}_{K_n}=1 \quad\Leftrightarrow\quad \Cl^{[\ell]}_{K_n}=1$ pour tout $n\in\NN$.}
\end{ThC}

Nous interrogeons ensuite plus particulièrement les conséquences de la trivialité du sous-module sauvage dans le cas abélien semi-simple en liaison avec les groupes d'unités et de classes circulaires introduits par Sinnott. La description explicite que nous obtenons (cf. {Théorème \ref{Thé}) fait apparaître une faille dans la preuve de la conjecture faible récemment proposée par Nguyen Quang Do dans \cite{Ng2} en correction de \cite{Ng1}.
\medskip

L'ensemble de ce travail peut ainsi être regardé comme une réinterprétation et un approfondissement des résultats classiques sur la conjecture de Greenberg en résonance avec des travaux récents de Gras \cite{Gra5} et de Nguyen Quang Do \cite{Ng3}. Par exemple, le cas semi-simple que nous traitons en fin de note donne un éclairage nouveau sur les travaux de Fukuda et Taya sur les corps quadratiques réels (cf. \cite{Fu1,Fu2,FuKo,FuTa}). De même, nombre de calculs de Gras s'interprètent de façon élégante en termes logarithmiques, notamment les calculs normiques qui font implicitement intervenir les valuations logarithmiques attachées aux places sauvages.

Enfin, toujours dans le cas complètement décomposé, nous montrons que la condition de trivialité du sous-module fini $\,\F_K^{\phantom{*}}$ de $\,\C_K^{\phantom{*}}$, qui se lit directement sur le sous-module sauvage et implique (au moins pour $\ell$-impair) la $\ell$-rationalité du corps (cf. Sco. \ref{rationnel}), permet de décomposer les $\ell$-adifiés des groupes de $\ell$-unités des divers étages de la tour cyclotomique (ainsi que leurs sous-groupes d'unités et, le cas échéant, d'unités circulaires) à l'aide de leurs sous-groupes logarithmiques et du seul groupe des unités de $K$.
\medskip

\noindent{\bf\em Remerciements.} Cette note a donné lieu à de nombreux échanges éclairants avec Karim Belabas, Jean-Robert Belliard, Georges Gras et Thong Nguyen Quang Do, que je remercie tout naturellement ici.

\newpage

%%%%%%%%%%%%%%%%%%%%%%%%%%%%%%%%%%%%%%%%%%%%%%%%%%%%%%%%%%%
\section{ Formulations de la conjecture de Greenberg}
%%%%%%%%%%%%%%%%%%%%%%%%%%%%%%%%%%%%%%%%%%%%%%%%%%%%%%%%%%%
\smallskip

Dans tout ce qui suit $\ell$  est un nombre premier arbitraire et  $K$ un corps de nombres totalement réel de degré $d$.\par

Nous sommes plus particulièrement intéressés au cas où la place $\ell$ se décompose complètement dans $K$. Néanmoins, certains résultats étant valables indépendamment de cette hypothèse, nous le précisons explicitement chaque fois qu'elle est utilisée.\medskip

\noindent$\bullet$ $K_\infty=\cup_{n\in\NN}K_n$ désigne la $\Zl$-extension cyclotomique de $K$ (avec la convention $[K_n:K]=\ell^n$); et $\Gamma=\gamma^{\Zl}=\Gal(K_\infty/K)$ son groupe de Galois; puis $\Lambda=\Zl[[\gamma-1]]$ l'algèbre d'Iwasawa associée.\smallskip

\noindent$\bullet$ $K_n^{nr}$ désigne la $\ell$-extension abélienne non-ramifiée $\infty$-décomposée maximale de $K_n$. Il suit donc:
\centerline{$\Gal(K_n^{nr}/K_n)\simeq\,\Cl_{K_n}$,}
\noindent où $\,\Cl_{K_n}$ est le $\ell$-sous-groupe de Sylow du groupe des classes d'idéaux (au sens ordinaire) de $K_n$.\smallskip

\noindent$\bullet$ $K_n^{\ell d}$ désigne la $\ell$-extension abélienne non-ramifiée $\ell\infty$-décomposée maximale de $K_n$. Il suit donc:
\centerline{$\Gal(K^{\ell d}/K_n)\simeq \,\Cl'_{K_n}$,}
\noindent où $\,\Cl'_{K_n}=\Cl_{K_n}/\Cl_{K_n}^{\,[\ell]}$ est le $\ell$-groupe des $\ell$-classes d'idéaux de $K_n$, i.e. le quotient du $\ell$-groupe $\,\Cl_{K_n}$ par le sous-groupe engendré par les classes des idéaux construits sur les premiers au-dessus de $\ell$.\smallskip

\noindent$\bullet$ $K_n^{lc}$ désigne la pro-$\ell$-extension abélienne localement cyclotomique maximale de $K_n$, i.e. la plus grande pro-$\ell$-extension abélienne de $K_n$ partout complètement décomposée sur $K_\infty$. Il suit donc:
\centerline{$\Gal(K^{lc}/K_\infty)\simeq \,\wCl_{K_n}$,}
\noindent où $\,\wCl_{K_n}$ est le $\ell$-groupe des classes logarithmiques (de degré nul) du corps $K_n$. La conjecture de Gross-Kuz'min (pour $K_n$ et $\ell$) postule précisément que ce dernier groupe est fini.\smallskip

\noindent$\bullet$ $K_n^{\ell r}$ désigne la pro-$\ell$-extension abélienne maximale de $K_n$ qui est $\ell$-ramifiée et $\infty$-décomposée, i.e. non-ramifiée en dehors des places au-dessus de $\ell$ et complètement décomposée aux places à l'infini. Le corps $K_n^{\ell r}$ contient naturellement la $\Zl$-extension cyclotomique $K_\infty$ et la conjecture de Leopoldt postule précisément qu'il est de degré fini sur $K_\infty$ puisque $K_n$ est totalement réel. Dans ce cas, nous écrivons $\T_{K_n}=\Gal(K_n^{\ell r}/K_\infty)$ le sous-groupe de torsion de $\Gal(K_n^{\ell r}/K_n)$. Il suit donc:
\centerline{$\Gal(K_n^{\ell r}/K_\infty)\simeq \,\T_{K_n}$.}\smallskip

\noindent$\bullet$ $K_n^{bp}$ désigne l'extension de Bertrandias-Payan associée à $K_n$, i.e. le compositum des $\ell$-extensions cycliques de $K_n$ qui sont {\em localement} $\Zl$-plongeables sur $K_n$. Le corps $K$ étant pris ici totalement réel, cette extension n'intervient dans la présente étude que pour $\ell=2$:
\begin{itemize}
\item Pour $\ell$ impair, $K_n^{bp}$ coïncide, en effet, avec $K^{\ell r}_n$ dès lors que les places $\l$ au-dessus de $\ell$ sont complètement décomposées dans $K$, puisque les complétés $K_{\l_n}$ des $K_n$ aux places $\l_n$ au-dessus de $\ell$ ne contiennent alors pas les racines $\ell$-ièmes de l'unité.
\item Pour $\ell=2$, supposé complètement décomposé dans $K$, la somme $\mu^{\phantom{*}}_{K_\ell}=\oplus_{\l|2}\, \mu^{\phantom{*}}_{K_\l}$ est un $\FF_2$-espace de dimension $d$ et il en est de même au niveau $K_n$, de sorte qu'il vient:

\centerline{$\Gal(K_n^{\ell r}/K_n^{bp}) \simeq \oplus_{\l_n|2}\, \mu^{\phantom{*}}_{K_{\l_n}}/ \mu^{\phantom{*}}_{K_n} \simeq \{\pm1\}^{d-1}$.}
\end{itemize}\smallskip

\noindent$\bullet$ $K_\infty^{nr}$ désigne la pro-$\ell$-extension abélienne non-ramifiée $\infty$-décomposée maximale de $K_\infty$. Il suit:
\centerline{$K_\infty^{nr}=\bigcup_{n\in\NN}K_n^{nr}$ \quad et \quad $\Gal(_\infty^{nr}/K_\infty) \simeq \,\C_K = \varprojlim \Cl_{K_n}$.}
Le groupe $\,\C_K$ est canoniquement un $\Lambda$-module noethérien et de torsion. On note $\lambda_K$ et $\mu_K$ ses invariants d'Iwasawa, i.e. le degré et la $\ell$-valuation de son polynôme caractéristique $\chi_{\C_K}\in\Zl[\gamma-1]$.\smallskip

\noindent$\bullet$ $K_\infty^{cd}$ dénote la pro-$\ell$-extension abélienne non-ramifiée et $\ell\infty$-décomposée maximale de $K_\infty$. Il suit:
\centerline{$K_\infty^{cd}=\bigcup_{n\in\NN}K_n^{\ell d}$ \quad et \quad $\Gal(K_\infty^{cd}/K_\infty)=\,\C'_K\simeq\varprojlim \Cl'_{K_n}$.}
Le groupe $\,\C'_K$, qui est appelé module de Kuz'min-Tate dans \cite{J55}, est le quotient de $\,\C_K$ par le sous-groupe $\,\C_K^{\,[\ell]}=\varprojlim \Cl_{K_n}^{\,[\ell]}$ construit sur les classes des premiers au-dessus de $\ell$. Et, comme la montée dans la $\Zl$-extension cyclotomique $K_\infty/K_n$ éteint toute possibilité d'inertie aux places étrangères à $\ell$ dans une $\ell$-extension de $K_\infty$, le corps $K_\infty^{cd}$ est aussi la plus grande pro-$\ell$-extension abélienne de $K_\infty$ qui est complètement décomposée partout.
\medskip

Enfin, nous utilisons librement dans l'ensemble de l'article les notations usuelles de la Théorie $\ell$-adique du corps de classes telle qu'exposée dans \cite{J28}. Eles sont introduite au fur et à mesure des besoins et récapitulées en appendice pour la commodité du lecteur.

\newpage

%%%%%%%%%%%%%%%%%%%%%%%%%%%%%%%%%%%%%%%%%%%%%%%%%%%%%%%%%%%
\subsection{Rappels sur le quotient de Herbrand}%\medskip
%%%%%%%%%%%%%%%%%%%%%%%%%%%%%%%%%%%%%%%%%%%%%%%%%%%%%%%%%%%

Si $\X$ est un $\Lambda$-module noethérien et de torsion dont le polynôme caractéristique $\chi\in\Zl [\gamma-1]$ n'est pas divisible par $(\gamma-1)$, le noyau $\X^\Gamma$ et le conoyau $^{\Gamma\!}\X$ de la multiplication par $\gamma-1$ dans $\X$ sont tous deux finis. On peut donc définir le quotient de Herbrand de $\X$ par:
$$
q(\X)=|\X^{\Gamma}|/|^{\Gamma\!}\X|
$$
Étant donnée une suite exacte courte de tels modules $1 \rightarrow \mathcal A \rightarrow \mathcal B \rightarrow \mathcal C \rightarrow 1$, le diagramme du serpent associée à la multiplication par $\gamma-1$ donne alors la suite exacte longue:\smallskip

\centerline{$1 \rightarrow \mathcal A^\Gamma \rightarrow \mathcal B^\Gamma \rightarrow \mathcal C^\Gamma \rightarrow 
 {}^{\Gamma\!}\mathcal A \rightarrow {}^\Gamma \mathcal B \rightarrow {}^\Gamma\mathcal C \rightarrow 1$,}\smallskip

\noindent qui fournit l'identité: $q(\mathcal B)=q(\mathcal A)q(\mathcal C)$. Or, si $\mathcal F$ est fini, la suite exacte canonique:

\centerline{$
\begin{CD}
1 @>>> \mathcal F^{\,\Gamma}  @>>> \mathcal F  @>\gamma-1>> \mathcal F  @>>>{}^\Gamma\mathcal F @>>> 1,
\end{CD}$}\smallskip

\noindent  donne immédiatement: $| \mathcal F^{\,\Gamma} | = |{}^\Gamma\mathcal F |$, i.e. $q(\mathcal F)=1$. Il suit de là que le quotient de Herbrand est invariant par pseudo-isomorphisme. On peut donc le calculer en remplaçant $\X$ par une somme directe de modules élémentaires $\Lambda/P\Lambda$.
Or, pour un tel module, on a banalement:\smallskip

\centerline{$\X^{\Gamma}=1 \qquad \rm{et} \qquad  {}^{\Gamma\!}\X = \Lambda/(P\Lambda+(\gamma-1)\Lambda)=\Lambda/\big(P(0)\Lambda+(\gamma-1)\Lambda\big)=\Zl/P(0)\Zl$.}\smallskip

\noindent De sorte que finalement il vient de façon générale: $q(\X)=(\Zl:\chi(0)\Zl)$. Et il suit:
$$
\X \sim 0 \Leftrightarrow \chi \in\Lambda^\times  \Leftrightarrow \chi(0)\in\Zl^\times  \Leftrightarrow q(\X)=1.
$$
\begin{Lem}\label{Herbrand}
Un $\Lambda$-module noethérien et de torsion dont le polynôme caractéristique n'est pas divisible par $\gamma-1$ est pseudo nul si et seulement si son quotient de Herbrand vaut 1.
\end{Lem}

Ce point acquis, intéressons-nous à la limite projective $\,\C_K=\varprojlim \Cl_{K_n}$ des $\ell$-groupes de classes.
Lorsque le corps de base $K$ est totalement réel, l'existence du quotient de Herbrand pour le groupe $\,\C_K$ est liée naturellement aux conjectures de Leopoldt et de Gross-Kuz'min (cf. e.g. \cite{J55,LMN,Ng3}):

\begin{Prop} \label{Conj}Soit $K$ un corps de nombres totalement réel. Alors:
\begin{itemize}\smallskip

\item[(i)] La conjecture de Leopoldt pour le corps $K$ (et le premier $\ell$) postule que la pro-$\ell$-extension abélienne $\ell$-ramifiée $\infty$-décomposée maximale $K^{\ell r}$ de $K$ est de degré fini sur $K_\infty$.\smallskip

\item[(ii)] Le polynôme caractéristique du $\Lambda$-module $\,\C_K=\varprojlim \Cl_{K_n}$ est étranger à $(\gamma-1)$ si et seulement si la pro-$\ell$-extension abélienne maximale $K^{ab}\cap K_\infty^{nr}$ de $K$ qui est non-ramifiée et $\infty$-décomposée sur $K_\infty$ est de degré fini sur $K_\infty$.\smallskip

\item[(iii)] La conjecture de Gross-Kuz'min pour le corps $K$ (et le premier $\ell$) postule que la pro-$\ell$-extension abélienne localement cyclotomique (i.e. complètement décomposée partout sur $K_\infty$) maximale $K^{lc}=K^{ab} \cap K_\infty^{cd}$ de $K$ est de degré fini sur $K_\infty$.\smallskip
\end{itemize}
Les inclusions immédiates $K^{lc} \subset K^{ab} \cap K_\infty^{nr} \subset K^{\ell r}$ hiérarchisent clairement ces trois conditions.
\end{Prop}

\begin{Sco}\label{Classes sauvages}
Lorsque les places au-dessus de $\ell$ ne se décomposent pas au-dessus de $K$ dans la tour cyclotomique $K_\infty/K$, la condition $(ii)$ impose au sous-groupe $\,\C_K^{\,[\ell]}=\varprojlim \Cl_{K_n}^{\,[\ell]}$ construit sur les classes des premiers au-dessus de $\ell$ d'être pseudo-nul; ce qui entraîne: $\Cl_{K_\infty}^{\,[\ell]}=\varinjlim \Cl_{K_n}^{\,[\ell]}=1$.\smallskip

La même conclusion vaut donc encore dès lors que le compositum $K_{d^{\phantom{*}}_0}$ des sous-corps de décomposition des places au-dessus de $\ell$ dans $K_\infty/K$ vérifie la conjecture de Leopoldt (en $\ell$); ce qui a lieu en particulier dès que $K$ est abélien sur $\QQ$.
\end{Sco}

\Preuve En l'absence de décomposition, le polynôme caractéristique $\chi_K^{[\ell]}$ du module $\,\C_K^{\,[\ell]}$ est une puissance de $(\gamma-1)$. La condition $(\gamma-1) \nmid \chi_K$  impose donc $\chi_K^{[\ell]}=1$; et $\,\C_K^{\,[\ell]}$ est alors pseudo-nul. \smallskip

Sous la conjecture de Leopoldt dans $K_{d_0}$ (en particulier lorsque $K$ est abélien sur $\QQ$) , nous pouvons appliquer la Proposition à l'extension $K_\infty/K_{d^{\phantom{*}}_0}$, ce qui nous donne encore $\,\C_K^{\,[\ell]}\sim 0$.

Dans tous les cas,  la norme $N_{m/n}$ est un isomorphime de $\, \Cl_{K_m}^{\,[\ell]}$ sur $\, \Cl_{K_n}^{\,[\ell]}$ pour $m \ge n \gg 0$; et l'identité $ N_{K_{m\!}/\!K_n}\circ j_{K_m/K_n}=\ell^{m-n}$ pour $m\ge n$ montre que $\,\Cl_{K_n}^{\,[\ell]}$ capitule dans $\,\Cl_{K_m}^{\,[\ell]}$ pour $m \gg  n$.

\newpage
%%%%%%%%%%%%%%%%%%%%%%%%%%%%%%%%%%%%%%%%%%%%%%%%%%%%%%%%%%%
\subsection{Formulations équivalentes classiques de la conjecture de Greenberg}
%%%%%%%%%%%%%%%%%%%%%%%%%%%%%%%%%%%%%%%%%%%%%%%%%%%%%%%%%%%

La conjecture de Greenberg est souvent énoncée pour les corps de nombres totalement réels qui satisfont la conjecture de Leopoldt pour un premier donné $\ell$. En fait, comme discuté dans la Proposition \ref{Conj}, il est tout à fait possible d'en donner des formulations indépendantes de cette dernière conjecture. C'est d'ailleurs ainsi que procède Greenberg dans \cite{Grb}.\smallskip

Rappelons que, pour chaque étage $K_n$ de la tour cyclotomique $K_\infty/K$, nous avons désigné par $\,\Cl'_{K_n}$ le $\ell$-groupe des $\ell$-classes d'idéaux de $K_n$, i.e. le quotient du $\ell$-groupe des classes $\Cl_{K_n}$ par le sous-groupe $\Cl_{K_n}^{\,[\ell]}$ engendré par les images des premiers au-dessus de $\ell$. Cela étant, nous avons:

\begin{Th}\label{CGreenberg}
Soient $K$ un corps de nombres totalement réel 
et $\,\C_K=\varprojlim \Cl_{K_n}$ la limite projective (pour la norme) des $\ell$-groupes de classes d'idéaux attachés aux étages finis $K_n$ de la $\Zl$-extension cyclotomique $K_\infty/K$. Si le polynôme caractéristique du $\Lambda$-module $\,\C_K$ n'est pas divisible par $\gamma-1$, les cinq assertions suivantes sont alors équivalentes et constituent la conjecture de Greenberg:\smallskip
\begin{itemize}
\item[(i)] Le quotient de Herbrand de $\,\C_K$ est trivial: $q(\C_K)=1$.\smallskip

\item[(ii)] Le $\Lambda$-module $\C_K$ est pseudo-nul: $\,\C_K \sim 0$.\smallskip

\item[(iii)] Les $\ell$-groupes de classes $\,\Cl_{K_n}$ sont d'ordres bornés (et donc stationnaires).\smallskip

\item[(iv)] Le $\ell$-groupe de classes $\,\Cl_{K_\infty}=\varinjlim \Cl_{K_n}$ est trivial:  $\,\Cl_{K_\infty}=1$.\smallskip

\item[(v)] Le sous-groupe ambige $\,\Cl_{K_\infty}^{\,\Gamma}=\varinjlim \Cl_{K_n}^{\,\Gamma}$ est trivial:  $\,\Cl_{K_\infty}^{\,\Gamma}=1$.\end{itemize}\smallskip

Lorsque, de plus, les places au-dessus de $\ell$ ne se décomposent pas dans $K_\infty/K$, ces mêmes conditions s'écrivent encore de façon équivalente:
\begin{itemize}\smallskip

\item[(vi)] Le quotient de Herbrand du module de Kuz'min-Tate $\,\C'_K=\varprojlim \Cl'_{K_n}$ est trivial: $q(\C'_K)=1$.\smallskip

\item[(vii)] Le  module de Kuz'min-Tate $\,\C'_K$ est pseudo-nul: $\,\C'_K \sim 0$.\smallskip

\item[(viii)] Les $\ell$-groupes de $\ell$-classes $\,\Cl'_{K_n}$ sont d'ordres bornés (et donc stationnaires).\smallskip

\item[(ix)] Le $\ell$-groupe des $\ell$-classes $\,\Cl'_{K_\infty}=\varinjlim \Cl'_{K_n}$ est trivial:  $\,\Cl'_{K_\infty}=1$.\smallskip

\item[(x)] Le sous-groupe ambige $\,\Cl'{\!}_{K_\infty}^{\;\Gamma}=\varinjlim \Cl'{\!}_{K_n}^{\;\Gamma}$ est trivial:  $\,\Cl'{\!}_{K_\infty}^{\;\Gamma}=1$.
\end{itemize}
\end{Th}

\Preuve Examinons successivement les cinq premières équivalences:\smallskip

$(i) \Leftrightarrow (ii)$ En vertu du lemme précédent.\smallskip

$(ii) \Leftrightarrow (iii)$ Puisque la norme $ N_{K_{m\!}/\!K_n}$ est ultimement surjective de $\,\Cl_{K_m}$ dans $\,\Cl_{K_m}$, dès que l'extension abélienne $K_m/K_n$ est totalement ramifiée en au moins une place.\smallskip

$(iii) \Leftrightarrow (iv)$ Si les ordres $|\Cl_{K_n}|$ sont bornés, les groupes de classes $\,\Cl_{K_n}$ sont ultimement isomorphes par la norme donc capitulent dans la tour puisque la composée de l'extension $j_{K_m/K_n}$ et de la norme $ N_{K_{m\!}/\!K_n}$ est l'exponentiation par le degré $[K_m:K_n]$. Et $\,\Cl_{K_\infty}$ est trivial. Inversement, si  $\,\Cl_{K_\infty}$ est trivial, les groupes de classes capitulent dans la tour, donc sont d'ordre borné, puisque la capitulation est structurellement bornée d'après Iwasawa. Et $\,\C_K$ est bien fini.

$(iv) \Leftrightarrow (v)$ Par action du pro-$\ell$-groupe $\Gamma$ sur la limite inductive de $\ell$-groupes $\,\Cl_{K_\infty}=\varinjlim \Cl_{K_n}$.\medskip

Ces équivalences acquises, la transposition aux $\ell$-groupes de $\ell$-classes $\,\Cl'_{K_n}$ ne pose aucune difficulté en vertu du Scolie \ref{Classes sauvages}, le sous-module $\,\C_K^{\,[\ell]}$ étant pseudo-nul par hypothèse.  On a ainsi: $(i)\Leftrightarrow (vi)$; $(ii)\Leftrightarrow (vii)$; $(vii)\Leftrightarrow (viii)$; et, comme plus haut, $(viii)\Leftrightarrow (ix)$, ainsi que $(ix)\Leftrightarrow (x)$.\medskip

\Remarque De façon générale, le polynôme caractéristique du $\Lambda$-module $\,\C_K=\varprojlim\Cl_{K_n}$ est le produit $\chi^{\phantom{l}}_K(\gamma-1)=\chi_K^{[\ell]}(\gamma-1)\chi_K'(\gamma-1)$ de ceux du sous-groupe $\,\C^{[\ell]}_K$ et du quotient $\,\C'_K=\C_K/\C_K^{[\ell]}$.
Lorsque les places au-dessus de $\ell$ ne se décomposent pas dans $K_\infty/K$, l'hypothèse $(\gamma-1) \nmid \chi^{\phantom{l}}_K(\gamma-1)$ entraîne la trivialité du premier facteur, en vertu du Scolie \ref{Classes sauvages}; de sorte que l'on a: $\chi^{\phantom{l}}_K= \chi_K'$.\medskip

Enfin, notons que la condition de non-décomposition des places au-dessus de $\ell$ dans $K_\infty/K$, n'est pas réellement restrictive: les places sauvages étant presque totalement ramifiées dans la $\Zl$-extension cyclotomique $K_\infty/K$, elle est toujours vérifiée dans $K_\infty/K_n$ pour $n\ge d_0^{\phantom{*}}$, si $K_ {d_0^{\phantom{*}}}$ désigne le compositum des sous-corps de décomposition des places au-dessus de $\ell$.

\newpage
%%%%%%%%%%%%%%%%%%%%%%%%%%%%%%%%%%%%%%%%%%%%%%%%%%%%%%%%%%%
\subsection{Introduction des $\ell$-groupes de classes logarithmiques}\medskip
%%%%%%%%%%%%%%%%%%%%%%%%%%%%%%%%%%%%%%%%%%%%%%%%%%%%%%%%%%%

Introduisons maintenant les $\ell$-groupes de classes logarithmiques $\,\wCl_{K_n}$ attachés aux divers étages de la tour $K_\infty/K$. Rappelons que ces groupes, qui sont décrits en termes de classes de diviseurs dans leurs corps de définition respectifs,  s'interprètent par la Théorie $\ell$-adique du corps de classes comme quotients des genres du pro-$\ell$-groupe $\,\C'_K = \varprojlim \,\Cl'_{K_n}$ (cf. e.g. \cite{J28,J31,J55}:\smallskip

\centerline{$\wCl_{K_n}\simeq\,\C'_K/\C'_K{}^{\omega_n}$, avec $\omega_n=\gamma^{\ell^n}-1$ pour tout $n\in\NN$.}\smallskip

Ainsi, du point de vue de Théorie d'Iwasawa, la question de la capitulation pour les $\ell$-groupes de classes logarithmiques se présente comme suit (cf. \cite{J55}, \S 8 et \cite{JM}, \S2.c): on dispose d'un module noethérien et de torsion $T$ sur l'algèbre d'Iwasawa $\Lambda = \Zl [[\gamma -1]]$ attachée à un groupe procyclique $\Gamma = \gamma^{\Zl}$~; on note $\omega_n =\gamma^{\ell^n} -1$; et on s'intéresse aux noyaux des morphismes de transition $j_{n,m}:T_n \rightarrow T_m^{\Gamma_{\!n}}$ pour $m \gg n$ induits par la multiplication par $\omega_m / \omega_n$ entre les quotients $T_n = T/\omega_n T$ et $T_m = T/\omega_m T$ pour $m-n$ assez grand; ce qui revient à considérer le noyau $\Ker_n = \Ker\,(T_n \rightarrow T_\infty)$ du morphisme d'extension $j_{n,\infty}$ à valeurs dans la limite inductive $T_\infty$ des $T_n$.\smallskip

On fait l'hypothèse (qui correspond à la conjecture de Gross-Kuz'min dans $K_\infty$) que le polynôme caractéristique $\chi(\gamma-1)$ du $\Lambda$-module $T$ est copremier à tous les $\omega_n$. \smallskip

Notant $F$ le plus grand sous-module fini de $T$ et $\ell^f$ son exposant, on obtient alors directement:\smallskip

\centerline{$\Ker_n = \{ t+\omega_n T \in T/ \omega_n T \, | \, (\omega_m / \omega_n)\ t  =
0 \} \underset{m \gg n}{=} (F + \omega_n T ) / \omega_n T$,}

\noindent puisque, pour $m > n$, les $\omega_m/\omega_n$ étant réputés copremiers avec le polynôme caractéristique de $T$, les éléments de $T$ qui sont tués par $\omega_m/\omega_n$ l'étant aussi par $\ell^f\chi(\gamma-1)$, ils sont forcément pseudo-nuls, i.e. contenus dans $F$; et, qu'inversement, les éléments de $F$ sont bien tués par  $\omega_m/\omega_n$ dès que $m-n$ est assez grand. En particulier, l'image $\bar T_n$ de $T_n$ dans $T_\infty=\varinjlim T_n$ est ainsi:\smallskip

\centerline{$\bar T_n \simeq T/(F+\omega_n T)$.}\smallskip 

Revenant alors au notations multiplicatives et appliquant ce qui précède au module de Kuz'min-Tate $\,\C'_K$,  notant enfin $\F'_K$ le sous-module fini de $\,\C'_K$, nous obtenons:\smallskip

\centerline{$j_{n,\infty}(\wCl_{K_n}) \simeq \,\C'_K/(\F'_K + \C'_K{}^{\omega_n})$.}

\begin{Th}\label{Finitude}
Soit $K$ un corps de nombres arbitraire et $\ell$ un nombre premier tel que la conjecture de Gross-Kuz'min soit vérifiée dans la $\Zl$-extension cyclotomique $K_\infty/K$ (i.e. pour lequel le polynôme caractéristique du module de Kuz'min-Tate $\,\C'_K=\varprojlim \,\Cl'_{K_n}$ est copremier avec les $\omega_n=\gamma^{\ell^n}-1$). Les propositions suivantes sont alors équivalentes:\smallskip

\begin{itemize}
\item[(i)] Le module de Kuz'min-Tate $\,\C'_K$ est fini : $\,\C'_K=\F'_K$.\smallskip

\item[(ii)] Le $\ell$-groupe des classes logarithmiques de $K_\infty$ est trivial: $\,\wCl_{K_\infty}=1$.\smallskip

\item[(iii)] Pour tout $n\in\NN$, le $\ell$-groupe des classes logarithmiques de $K_n$ capitule dans $K_\infty$.\smallskip

\item[(iv)] Le $\ell$-groupe des classes logarithmiques de $K$ capitule dans la tour $K_\infty/K$.
\end{itemize}
\end{Th}

\Preuve Observons que le module de Kuz'min-Tate $\,\C'_K$ est encore la limite projective $\varprojlim\,\wCl_{K_n}$.\par

$(i)$ S'il est fini, on a: $\,\C'_K\simeq\,\wCl_{K_n}$, pour $ n \gg 0$ et l'identité $ N_{K_{m\!}/\!K_n}\circ j_{K_m/K_n}(\wCl_{K_n})=\wCl_{K_n}^{\ell^{m-n}}\!=1$ pour $m\gg n$ montre que $\,\wCl_{K_n}$ capitule dans $\,\wCl_{K_m}$ pour $m \gg  n$.  Il suit: $\,\wCl_{K_\infty}=1$.\par

$(ii)$ Si $\,\wCl_{K_\infty}$ est trivial, tous les $\,\wCl_{K_n}$ (et, en particulier, $\;\wCl_K$) capitulent dans $K_\infty$.\par

$(iv)$ Si $\,\wCl_K$ capitule dans $K_\infty$, on a: $\,\C'_K/(\F'_K + \C'_K{}^{\,\omega})=1$; donc, par Nakayama: $\,\C'_K=\F'_K$.

\begin{Sco}
Sous les hypothèses du Théorème, pour tout $n \in \NN$, l'image dans $\,\wCl_{K_\infty}$ du $\ell$-groupe des classes logarithmiques du sous-corps $K_n$ de $K_\infty$ fixé par $\,\Gamma_{\! n}$ est le sous-module des points fixes:\smallskip

\centerline{$j_{n,\infty}(\wCl_{K_n}) = \,\wCl_{K_\infty}^{\,\Gamma_{\!n}}$.}
\end{Sco}

\Preuve Revenons pour simplifier aux notations additives et notons $\bar T= T/F$ le quotient de $T$ par son sous-module fini $F$. D'après ce qui précède, l'image $\bar T_n$ de $T_n$ dans $T_\infty$ est alors:\smallskip

\centerline{$\bar T_n \simeq T/(F+\omega_n T) \simeq \bar T/\omega_n\bar T \simeq \frac{1}{\omega_n}\bar T/\bar T$;}

\noindent et il suit, comme annoncé:

\centerline{$T_\infty= \big(\bigcup_{n\in\NN}  \frac{1}{\omega_n}\bar T\big)/\bar T$ \quad et \quad $\bar T_n= T_\infty^{\Gamma_{\! n}}$.}

\newpage
%%%%%%%%%%%%%%%%%%%%%%%%%%%%%%%%%%%%%%%%%%%%%%%%%%%%%%%%%%%
\subsection{Interprétation logarithmique de la conjecture de Greenberg}\medskip
%%%%%%%%%%%%%%%%%%%%%%%%%%%%%%%%%%%%%%%%%%%%%%%%%%%%%%%%%%%

Le Théorème \ref{Finitude} fournit immédiatement le critère logarithmique suivant, indépendant de toute hypothèse sur la décomposition des places sauvages dans le corps totalement réel considéré, et qui peut être regardé comme la généralisation naturelle du critère de capitulation obtenu par Greenberg en présence d'une unique place au-dessus de $\ell$ (cf. \cite{Grb}, Th. 1):

\begin{Th}[Critère logarithmique]\label{CGLog}
Soient $K$ un corps de nombres totalement réel, $\ell$ un nombre premier et $K_\infty$ la $\Zl$-extension cyclotomique de $K$. Sous la conjecture de Leopoldt dans $K_\infty$ (donc en particulier dès que $K$ est abélien sur $\QQ$), le corps $K$ vérifie la conjecture de Greenberg (pour le premier $\ell$) si et seulement si le $\ell$-groupe des classes logarithmiques de $K$ capitule dans $K_\infty$.
\end{Th}

Comme indiqué dans l'introduction, ce résultat ramène la conjecture de Greenberg à une pure question de capitulation pour les classes logarithmiques. Outre son intérêt théorique, il ouvre de ce fait sur la possibilité de vérifier numériquement cette conjecture dans tous les cas où le contexte arithmétique permet de mettre en \oe uvre les algorithmes présentés dans \cite{BJ} et \cite{DJ+}.
Signalons enfin qu'une preuve indépendante vient d'en être donnée par Nguyen Quang Do (cf. \cite{Ng3}, Th. 2.1).\smallskip

\Preuve Les divers étages $K_n$ de la $\Zl$-tour $K_\infty/K$ sont encore totalement réels et vérifient ainsi la conjecture de Gross-Kuz'min (pour le premier $\ell$), puisqu'ils vérifient par hypothèse celle de Leopoldt. La capitulation du $\ell$-groupe $\,\wCl_K$ caractérise donc la finitude de $\,\C'_K$, qui équivaut à celle de $\,\C_K$, puisque les sous-groupes sauvages $\,\Cl_{K_n}^{[\ell]}$ sont alors bornés d'après le Scolie \ref{Classes sauvages}.
\medskip

Récapitulant tout ce qui précède, nous avons ainsi (en parallèle avec \cite{Ng3}, Th. 2.1):

\begin{Sco}\label{ThPrincipal}
Soit $K$ un corps de nombres totalement réel et $\ell$ un nombre premier.
\begin{itemize}\smallskip

\item[(i)] Si la $\Zl$-extension cyclotomique $K_\infty/K$ satisfait la conjecture de Gross-Kuz'min pour $\ell$, i.e. si le polynôme caractéristique du $\Lambda$-module de Kuz'min-Tate $\,\C'_K=\varprojlim\,\Cl'_{K_n}$ est copremier avec les $\omega_n=\gamma^{\ell^n}-1$ (et donc en particulier si $K$ est abélien sur $\QQ$), on a les équivalences:\smallskip
\begin{itemize}
\item[(i,a)]  $\,\C'_K$ est fini si et seulement si le $\ell$-groupe logarithmique $\,\wCl_K$ capitule dans $K_\infty$.

\item[(i,b)]  $\,\C'_K$ est trivial si et seulement si le $\ell$-groupe logarithmique $\,\wCl_K$ est trivial.
\end{itemize}\smallskip

\item[(ii)] Si  le compositum  $K_{d_0^{\phantom{*}}}$ des sous-corps de décomposition des places au-dessus de $\ell$ dans $K_\infty/K$ satisfait la conjecture de Leopoldt  (et donc en particulier si $K$ est abélien sur $\QQ$), le groupe $\,\C_K=\varprojlim\,\Cl_{K_n}$ et le module de Kuz'min-Tate $\,\C'_K$ sont simultanément finis ou pas.\smallskip

\item[(iii)] Lorsque ces deux conditions sont réunies et donc, en particulier, sous la conjecture de Leopoldt pour $\ell$ dans $K_\infty$, le corps $K$ satisfait la conjecture de Greenberg pour le premier $\ell$ si et seulement si le $\ell$-groupe des classes logarithmiques $\,\wCl_K$ capitule dans la tour $K_\infty/K$.
\end{itemize}
\end{Sco}

\begin{Sco}[Premier critère de Greenberg]
Soit $K$ un corps de nombres totalement réel qui admet une unique place au-dessus de $\ell$. Si celle-ci est totalement ramifiée dans $K_\infty/K$, le corps $K$ satisfait la conjecture de Greenberg (pour le premier $\ell$) si et seulement si le $\ell$-groupe des classes d'idéaux $\,\Cl_K$ capitule dans $K_\infty$.
\end{Sco}

\Preuve  La $\Zl$-tour cyclotomique $K_\infty/K$ étant supposée ici totalement ramifiée, le morphisme naturel de $\,\wCl^{\phantom{*}}_{K_n}$ vers $\,\Cl'_{K_n}$ est surjectif pour tout $n\in\NN$. Comme, de plus, il n'existe par hypothèse qu'une seule place $\l_n$ au-dessus de $\ell$ à chaque étage fini $K_n$ de la tour, chacun des corps $K_n$ satisfait la conjecture de Kuz'min-Tate (pour le premier $\ell$) et l'on a même, plus précisément:\smallskip

\centerline{$\wCl^{\phantom{*}}_{K_n} = \,\Cl'_{K_n}$\quad pour tout $n\in\NN$.}\smallskip

\noindent D'après le Théorème \ref{Finitude}, la capitulation du $\ell$-groupe logarithmique $\,\wCl_K$ assure la finitude du module de Gross-Kuz'min $\,\C'_K$. Tout revient donc à vérifier la finitude du sous-groupe sauvage $\,\C_K^{[\ell]}$ de $\,\C^{\phantom{*}}_K$ indépendamment de la conjecture de Leopoldt; en d'autres termes, que les $\,\Cl^{[\ell]}_{K_n}$ restent bornés dans la tour. Or, ces groupes sont ambiges puisque engendrés par les classes d'idéaux ambiges, et la formule des $\ell$-classes ambiges de Chevalley appliquée dans l'extension cyclique $K_n/K$ nous donne ici:\smallskip

\centerline{$|\,\Cl_{K_n}^{\,\Gamma}|\,=\,|\,\Cl^{\phantom{*}}_K|\,\frac{e_\l(K_n/K)}{[K_n:K]\big(E_K:E_K\cap N_{K_{\scale{.6}{\scriptstyle n}}/K}(K_n^\times)\big)}\,=\,|\,\Cl^{\phantom{*}}_K|$.}

%%%%%%%%%%%%%%%%%%%%%%%%%%%%%%%%%%%%%%%%%%%%%%%%%%%%%%%%%%%
%%%%%%%%%%%%%%%%%%%%%%%%%%%%%%%%%%%%%%%%%%%%%%%%%%%%%%%%%%%
\section{Étude du cas complètement décomposé}\smallskip
%%%%%%%%%%%%%%%%%%%%%%%%%%%%%%%%%%%%%%%%%%%%%%%%%%%%%%%%%%%
%%%%%%%%%%%%%%%%%%%%%%%%%%%%%%%%%%%%%%%%%%%%%%%%%%%%%%%%%%%

Avant d'aller plus loin sur la question qui nous préoccupe, précisons quelques éléments de la Théorie $\ell$-adique du corps de classes lorsque les places sauvages se décomposent complètement dans l'extension $K/\QQ$ que nous supposerons totalement réelle. Les notations sont celles de \cite{J31}.\medskip

Notons $\R_{K_\p}=\varprojlim \,K^\times/K^{\times\ell^n}$ le $\ell$-adifié du groupe multiplicatif du complété $K_\p$ de $K$ en $\p$; désignons par $\,\U^{\phantom{*}}_{K_\p}$ son sous-groupe unité (au sens habituel); et notons $\,\wU^{\phantom{*}}_{K_\p}$ le groupe des unités logarithmiques, i.e. le  sous-groupe de normes attaché à la $\Zl$-extension cyclotomique $K_\p^c$ de $K_\p$: \smallskip
\begin{itemize}

\item Pour $\p$ {\em modérée}, i.e. pour $\p\nmid\ell\infty$, le groupe  $\R_{K_\p}$ est le produit direct  $\R_{K_\p}= \mu_{K_\p}\pi_\p^{\Zl}$ du $\ell$-groupe des racines de l'unité contenues dans $K_\p$ et du $\Zl$-module libre construit sur l'image $\pi_\p$ d'une uniformisante arbitraire.
Ses sous-groupes unités sont ainsi: $\,\U^{\phantom{*}}_{K_\p}=\,\wU^{\phantom{*}}_{K_\p}=\mu^{\phantom{*}}_{K_\p}$.\smallskip

\item Pour $\l$ {\em sauvage}, i.e. pour $\l\mid\ell$, c'est le produit $\R_{K_\l}=\,\U_{K_\l}\pi_\l^{\Zl}$ du groupe des unités principales du complété $K_\l$ et du $\Zl$-module libre construit sur une uniformisante $\pi_\l$. Si $\ell$ se décompose complètement dans $K$, il est naturel de prendre pour $\pi_\l$ l'image de $\ell$ dans $K_\l$. De plus:
\begin{itemize}
\item si $\ell$ est impair, le groupe $\,\U_{K_\l}$ des unités principales est lui-même un $\Zl$-module libre: il est engendré par l'image $\wi\pi_\l$ de $1+\ell$; et il vient ainsi: $\R_{K_\l}=\,\U_{K_\l}\,\wU_{K_\l}=\wi\pi_\l^{\Zl}\pi_\l^{\Zl}$;
\item si $\ell$ vaut 2, il vient: $\R_{K_\l}=\{\pm 1\}\,\wi\pi_\l^{\Zl}\,\pi_\l^{\Zl}$; $\U_{K_\l}=\{\pm 1\}\,\wi\pi_\l^{\Zl}$; et $\,\wU_{K_\l}=\{\pm 1\}\,\pi_\l^{\Zl}$.
\end{itemize}\smallskip

\item Pour $\p$ réelle enfin, on a: $\R_{K_\p}=\mu_{K_\p}=\{\pm 1\}$, si $\ell$ vaut 2;  et $\R_{K_\p}=\mu_{K_\p}=1$, sinon.
\end{itemize}\smallskip

Le $\ell$-adifié du groupe des idèles du corps $K$ est le produit restreint $\J_K=\prod_\p^{res}\R_{K_\p}$ formé des familles $(x_\p)_\p$ dont presque tous les éléments sont des unités. Son sous-groupe unité est le produit $\,\U_K=\prod_\p\U_{K_\p}$; le sous-groupe des unités logarithmiques de $\J_K$ est le produit $\,\wU_K=\prod_\p\wU_{K_\p}$. \par

Nous notons  $\R_{K_\ell}=\prod_{\l\mid\ell}\R_{K_\l}$ puis $\,\U_{K_\ell}=\prod_{\l\mid\ell}\U_{K_\l}$ et  $\,\wU_{K_\ell}=\prod_{\l\mid\ell}\wU_{K_\l}$ les groupes semi-locaux.\medskip

Le sous-groupe des idèles principaux est le tensorisé $\R_K=\Zl\otimes_\ZZ K^\times$ du groupe multiplicatif de $K$; il se plonge canoniquement dans $\J_K$ via le plongement diagonal de $K$ dans le produit de ses complétés.
Le quotient $\J_K/\R_K$ est un groupe compact et la Théorie $\ell$-adique du corps de classes l'identifie au groupe de Galois $\Gal(K^{ab}/K)$ de la pro-$\ell$-extension abélienne maximale de $K$. La Théorie de Galois établit alors une bijection entre les sous-extensions de $K^{ab}/K$ et les sous-groupes fermés de $\J_K$ qui contiennent $\R_K$.
Le sous-groupe normique attaché à la $\Zl$-extension cyclotomique $K_\infty$ de $K$ est ainsi le noyau $\wJ_K$ de la formule du produit pour les valeurs absolues.\smallskip

Le sous-groupe unité $\E_K=\R_K\cap\,\U_K$ de $\R_K$ est le tensorisé $\ell$-adique $\E_K=\Zl\otimes_\ZZ E_K$ du groupe des unités de $K$. Le groupe des unités logarithmiques $\wE_K=\R_K\cap\,\wU_K$ est le sous-groupe des normes cyclotomiques dans $\R_K$; contrairement à $\E_K$, ce n'est pas en général le $\ell$-adifié d'un sous-groupe de $K^\times$. Mais $\wE_K$ et $\E_K$ sont des sous-groupes du $\ell$-adifié $\E'_K=\Zl\otimes_\ZZ E'_K$ du groupe $E'_K$ des $\ell$-unités de $K$. Et la conjecture de Leopoldt (pour $K$ et pour $\ell$) affirme l'injectivité du morphisme de semi-localisation $s^{\phantom{*}}_\ell$ de $\E'_K$ dans $\R_{K_\ell}$ induit par le plongement de $K^\times$ dans $K_\ell^\times=\prod_{\l\mid\ell}K_\l^\times$.\smallskip

Lorsque, en outre, les places au-dessus de $\ell$ se décomposent complètement dans $K/\QQ$, l'expression explicite des valeurs absolues montre que l'intersection de $\wJ_K$ avec $\,\U_{K_\ell}=\prod_{\l\mid\ell}\U_{K_\l}$ coïncide avec la pré-image $\,\U^*_{K_\ell}$ dans  $\,\U_{K_\ell}$ du sous-groupe de torsion $\mu^{\phantom{*}}_{\Ql}\simeq\Zl/2\Zl$ par la norme $N_{K/\QQ}$:

\centerline{$\,\U^*_{K_\ell}=\{(x_\l)_\l\in\U_{K_\ell} \; |\; \prod_{\l\mid\ell}x_\l \in\mu^{\phantom{*}}_{\Ql}\}$.}

\noindent Et la conjecture de Leopoldt affirme alors que l'image $\E_{K_\ell}{=}s^{\phantom{*}}_\ell(\E_K)$ de $\E_K$ dans $\,\U_{K_\ell}$ est d'indice fini dans $\,\U^*_{K_\ell}$ ou, de façon équivalente, que  l'image $\E'_{K_\ell}{=}s^{\phantom{*}}_\ell(\E'_K)$ de $\E'_K$ dans $\R_{K_\ell}$ est d'indice fini dans le produit $\R^*_{K_\ell}=\U^*_{K_\ell}\,\wU_{K_\ell}$. Comme expliqué plus haut, elle entraîne dans ce cas la conjecture de Gross-Kuz'min, dont une formulation revient à postuler de même que que l'image $\wE_{K_\ell}{=}s^{\phantom{*}}_\ell(\wE_K)$ de $\wE_K$ dans $\R_{K_\ell}$ est d'indice fini dans $\,\wU_{K_\ell}$ (cf. \cite{J55}).\smallskip

Enfin, du fait que les places au-dessus de $\ell$ sont totalement ramifiées dans $K_\infty/K$, chaque classe dans le $\ell$-groupe des classes d'idéaux au sens ordinaire peut être représentée par un idéal de degré nul (ou, si l'on préfère par un idèle de $\wJ_K$) de sorte qu'on a: $\,\Cl_K=\J_K/\U_K\R_K \simeq \wJ_K/\U^*_K\R_K$, analoguement au $\ell$-groupe des classes logarithmiques défini par l'identité: $\,\wCl_K=\wJ_K/\wU_K\R_K$. Tous deux admettent comme quotient le $\ell$-groupe des $\ell$-classes donné par: $\,\Cl'_K=\wJ_K/\R^*_{K_\ell}\prod\mu_\p{\phantom{*}}\R_K$.

\newpage
%%%%%%%%%%%%%%%%%%%%%%%%%%%%%%%%%%%%%%%%%%%%%%%%%%%%%%%%%%%
\subsection{Calcul du quotient des genres dans le cas $\ell$-décomposé}\label{2.1}
%%%%%%%%%%%%%%%%%%%%%%%%%%%%%%%%%%%%%%%%%%%%%%%%%%%%%%%%%%%

Supposons désormais que les places au-dessus de $\ell$ se décomposent complètement dans $K$.\smallskip

Notons $K^{bp}$ l'extension de Bertrandias-Payan du corps $K$, i.e. le compositum des $\ell$-extensions cycliques de $K$ qui sont {\em localement} plongeables dans une $\Zl$-extension. Par la Théorie $\ell$-adique du corps de classes, l'extension $K^{bp}$ est associée au groupe de normes $\prod_\p\,\mu_{K_\p}\,\R_K$ (cf. e.g. \cite{J31}, Ex. 2.9). En particulier, elle est $\ell$-ramifiée (i.e. non-ramifiée en dehors de $\ell$) et $\infty$-décomposée (i.e. totalement réelle), puisque fixée par $\,\U_{K_\p}=\mu_{K_\p}$ pour $\p$ modérée et par $\R_{K_\p}=\mu_{K_\p}$ pour $\p$ réelle.\par
Aux places $\l | \ell$, le sous-groupe d'inertie $I_\l(K^{bp}/K)$ s'identifie à l'image de $\,\U_{K_\l}$ dans le quotient $\J_K/\prod_\p\,\mu_{K_\p}\,\R_K$.
Or, nous avons ici: $\,\U_{K_\l}=\mu_{K_\l}\wi\pi_\l^{\Zl}$; et $I_\l(K^{bp}/K)$ est procyclique, engendré par l'image de $\wi\pi_\l$, qui est d'ordre infini puisque la sous-extension cyclotomique $K_\infty/K$ est infiniment ramifiée en $\l$.
Et $K^{bp}/K_\infty$ est non-ramifiée aux places $\l\mid\ell$.

Généralisant un résultat classique d'Ozaki et Taya (cf. \cite{Oz,OT,Ta}), nous avons ainsi:

\begin{Th}\label{Genres}
Soient $K$ un corps de nombres totalement réel, $\ell$ un nombre premier complètement décomposé dans $K$, et $K^{bp}$ la pro-$\ell$--extension de Bertrandias-Payan attchée à $K$. Alors: 
\begin{itemize}
\item Pour $\ell\ne 2$, la pro-$\ell$-extension abélienne $\ell$-ramifiée $\infty$-décomposée maximale $K^{\ell r}$ coïncide avec $K^{bp}$ et c'est la plus grande pro-$\ell$-extension abélienne $K$ qui est non-ramifiée sur $K_\infty$.
\item Pour $\ell=2$ et sous la conjecture de Leopoldt, on a $[K^{\ell r}:K^{bp}]=2^{d-1}$ et la plus grande pro-2-extension abélienne totalement réelle de $K$ non-ramifiée sur $K_\infty$ est encore $K^{bp}$.
\end{itemize}
Dans les deux cas,  le groupe $\,\T_K^{bp}=\Gal(K^{bp}/K_\infty)$ s'identifie au quotient des genres $\,^\Gamma\C_K$ attaché à $\,\C_K$ et la conjecture de Leopoldt (pour le corps $K$ et le premier $\ell$) postule donc précisément la finitude de $\,^\Gamma\C_K$: en particulier, les conditions $(i)$ et $(ii)$ de la Proposition \ref{Conj} sont alors équivalentes.
\end{Th}

\Preuve On a, en effet: $\Gal(K^{\ell r}/K^{bp})\simeq \,\mu^{\phantom{*}}_{K_\ell}/s_\ell^{\phantom{*}}(\mu_K^{\phantom{*}})      \underset{^\text{déf}}{=}\big(\prod_{\l\mid\ell}\,\mu^{\phantom{*}}_{K_\l}\big)/s_\ell^{\phantom{*}}(\mu^{\phantom{*}}_K)$  (cf.  \cite{J31}, Ex. 2.9).\medskip

Considérons alors le schéma de corps:

\begin{center}
\unitlength=1.5cm
\begin{picture}(6.6,2.8)

\put(0.7,0){$K$}
\put(0.8,0.3){\line(0,1){1.5}}
\put(0.6,2){$K_\infty$}
\put(3.2,0){$K^{nr}$}
\put(1.2,0.05){\line(1,0){1.6}}

\bezier{60}(0.6,0.3)(0.4,1.2)(0.6,1.8)
\put(0.3,1.0){$\Gamma$}

\put(1.2,2.05){\line(1,0){1.5}}
\put(2.9,2){$K^{nr}K_\infty$}
\put(3.9,2.05){\line(1,0){1.5}}
\put(5.6,2){$K^{bp}$}
\put(3.3,0.3){\line(0,1){1.5}}

\bezier{100}(1.2,2.3)(3,2.7)(5.4,2.3)
\put(3.2,2.65){$\T_K^{bp}$}

\bezier{50}(1.2,1.8)(2.0,1.6)(2.9,1.8)
\put(1.9,1.4){$\Cl_K$}

\bezier{50}(3.7,1.8)(4.5,1.6)(5.4,1.8)
\put(3.9,1.4){$\U^*_{K_\ell}/s^{\phantom{*}}_\ell(\E^{\phantom{*}}_K)\mu^{\phantom{*}}_{K_\ell}$}

\end{picture}
\end{center}

La Théorie $\ell$-adique du corps de classes (cf. \cite{J31}, \S2) nous donne l'isomorphisme:\smallskip

\centerline{$\Gal(K^{bp}/K^{nr}) \simeq \,\U^{\phantom{*}}_{K_\ell}/s^{\phantom{*}}_\ell(\E^{\phantom{*}}_K)\mu^{\phantom{*}}_{K_\ell}$;
\quad puis \quad $\Gal(K^{bp}/K^{nr}K_\infty) \simeq \,\U^*_{K_\ell}/s^{\phantom{*}}_\ell(\E^{\phantom{*}}_K)\mu^{\phantom{*}}_{K_\ell}$}\smallskip

\noindent  avec $\,\U^*_{K_\ell}= \,\U^{\phantom{*}}_{K_\ell}\cap\wJ_K$ (cf. supra). 
Notant alors $\ell^e$ l'exposant du quotient $ \,\U^*_{K_\ell}/s^{\phantom{*}}_\ell(\E^{\phantom{*}}_K)\mu^{\phantom{*}}_{K_\ell}$, nous obtenons, pour tout $n \ge e$, l'isomorphisme:\smallskip

\centerline{$\U^*_{K_\ell}/s^{\phantom{*}}_\ell(\E^{\phantom{*}}_K)\mu^{\phantom{*}}_{K_\ell} \simeq (s^{\phantom{*}}_\ell(\E^{\phantom{*}}_K)\cap\,\U_{K_\ell}^{\ell^n})/s^{\phantom{*}}_\ell(\E_K^{\ell^n})$.}\smallskip

\noindent Pour interpréter le numérateur, observons maintenant que $\,\U_{K_\ell}^{\ell^n}$ est précisément le sous groupe normique de $\,\U_{K_\ell}$ relatif à l'extension $K_n/K$; de sorte que les éléments de l'intersection avec $s_\ell(\E_K)$ sont exactement les images des unités de $K$ qui sont normes (locales comme globales) dans l'extension cyclique $K_n/K$. Le morphisme $s_\ell$ étant injectif sous la conjecture de Leopoldt, il suit:

\begin{Cor}\label{CGenres}
Lorsque les places au-dessus de $\ell$ se décomposent complètement dans $K$ et sous la conjecture de Leopoldt dans $K$, l'ordre du groupe de torsion $\,\T^{bp}_K$ est donné, pour tout $n\ge e$,  par:\smallskip

\centerline{$|\T^{bp}_K| = |\,\Cl_K| \big((\E_K\cap N_{K_{n\!}/\!K}(\R_{K_n}):\E_K^{\ell^n}\big) =  |\,\Cl_K|\;\frac{\ell^{n(d-1)}}{(\E_K:\,\E_K\cap\,   N_{K_{n\!}/\!K}(\R_{K_n}))}$}\smallskip

\noindent C'est la valeur du nombre de genres $|{}^\Gamma \Cl_{K_n}|$ dans l'extension cyclique $K_n/K$; de sorte que $K^{bp}$ est encore la plus grande $\ell$-extension abélienne de $K$ qui est non-ramifiée et $\infty$-décomposée sur $K_n$.
\end{Cor}

\newpage
%%%%%%%%%%%%%%%%%%%%%%%%%%%%%%%%%%%%%%%%%%%%%%%%%%%%%%%%%%%
\subsection{Critère de Greenberg dans le cas $\ell$-décomposé}
%%%%%%%%%%%%%%%%%%%%%%%%%%%%%%%%%%%%%%%%%%%%%%%%%%%%%%%%%%%

Revenant alors aux sous-groupes ambiges, nous obtenons immédiatement le résultat suivant:

\begin{Lem}
Soit $K$ un corps de nombres totalement réel qui vérifie la conjecture de Leopoldt pour le premier $\ell$ et où les places au-dessus de $\ell$ se décomposent complètement. Alors $K$ satisfait les conditions équivalentes de la conjecture de Greenberg listées plus haut si et seulement si la norme arithmétique $ N_{K_{m\!}/\!K_n}:\;\Cl_{K_m}^{\,\Gamma}\rightarrow\Cl_{K_n}^{\,\Gamma}$ réalise un isomorphisme pour tous $m\ge n\gg 0$.
\end{Lem}

\Preuve D'après le Théorème \ref{Genres} ci-dessus, les quotients des genres  ${}^\Gamma \Cl_{K_n}$ ont ultimement le même ordre $|\T^{bp}_K| = |^\Gamma\C_K|$ dans la tour $K_\infty/K$. Il en va donc de même des sous-groupes ambiges $\,\Cl_{K_n}^{\,\Gamma}$, en vertu de l'égalité: $|\,\Cl_{K_n}^{\,\Gamma}|=|{}^\Gamma\Cl_{K_n}|$, dans l'extension cyclique $K_n/K$. L'égalité $|^\Gamma\C_K|= |\,\C_K^{\,\Gamma}|$ a donc lieu si et seulement si les sous-groupes ambiges $\,\Cl_{K_n}^{\,\Gamma}$ s'identifient ultimement à leur limite projective, i.e. lorsqu'ils sont ultimement isomorphes par la norme. Il suit:

\begin{Th}[Critère de Greenberg]\label{CrGreenberg}
Sous la conjecture de Leopoldt pour $K$, dès lors que les places au-dessus de $\ell$ se décomposent complètement dans $K/\QQ$, les dix conditions équivalentes du Théorème \ref{CGreenberg} sont réalisées si et seulement si les conditions équivalentes suivantes le sont également:\smallskip

\begin{itemize}
\item[(xi)] Les sous-groupes ambiges $\,\Cl_{K_n}^{\,\Gamma}$ sont ultimement isomorphes: $\,\Cl_{K_n}^{\,\Gamma}= N_{K_{m\!}/\!K_n}(\Cl_{K_n}^{\,\Gamma})$.\smallskip

\item[(xii)] Pour tout $n\gg 0$, le sous-groupe ambige $\,\Cl_{K_n}^{\,\Gamma}$ de $\,\Cl_{K_n}$ est engendré par les classes des premiers au-dessus de $\ell$:  $\,\Cl_{K_n}^{\,\Gamma}=\Cl_{K_n}^{\,[\ell]}$.

\item[(xiii)] Pour tout $n\gg 0$, on a les deux conditions: $\displaystyle
\begin{cases}
\text{(xiii,a)}\quad \Cl'_K \text{ capitule dans } \,\Cl'_{K_n},\smallskip\\
\text{(xiii,b)}\quad E_K\cap N_{K_{n\!}/\!K}(K_n^\times)= N_{K_{n\!}/\!K}(E_{K_n}),
\end{cases}$
\end{itemize}
\end{Th}

\Preuve La condition $(xii)$ constitue à proprement parler le critère donné par Greenberg.
\smallskip

 L'équivalence avec $(xiii)$ provient, d'une part, de  l'isomorphisme classique de Chevalley comparant classes ambiges et classes d'idéaux ambiges, unités qui sont normes et normes d'unités: 
 
 \centerline{$\,\Cl_{K_n}^{\,\Gamma}/\Cl^{\phantom{l}}_{K_n}(Id^{\,\Gamma}_{K_n})\simeq \big( E_K\cap N_{K_{n\!}/\!K}(K_n^\times)\big) /N_{K_n/K}(E_{K_n})$;}
 \smallskip
 
\noindent et, d'autre part, de l'égalité $\,\Cl_{Kn}^{\phantom{l}}(Id^{\,\Gamma}_{K_n})=\Cl_{K_n}^{\,[\ell]}$, qui traduit la capitulation de $\,\Cl'_K$ dans $\Cl'_{K_n}$.\smallskip
 
Compte tenu du Lemme, seule reste donc à vérifier l'équivalence avec la condition $(xi)$. Or:
\begin{itemize}

\item d'un côté, l'égalité immédiate $\,\Cl_{K_n}^{[\ell]}= N_{K_{m\!}/\!K_n}(\Cl_{K_n}^{[\ell]})$, pour $m\ge n$, les places au-dessus de $\ell$ étant totalement ramifiées dans $K_\infty/K$, nous donne l'implication: $(xii) \Rightarrow (xi)$;\smallskip

\item d'un autre côté, sous la condition $(xi)$, le groupe $\,\Cl_{K_\infty}$ est trivial d'après le Lemme et les idéaux de $K$ se principalisent dans $K_n$ pour $n \gg 0$. Prenons aussi $n$ assez grand pour avoir l'isomorphisme $(xi)$; choisissons $m=n+e$, où $e$ est l'indice défini dans la section précédente; et  partons d'un idéal $\a_n= N_{K_{m\!}/\!K_n}(\a_m)$ représentant une classe de $\,\Cl_{K_n}^{\,\Gamma}$ avec $\a_m$ appartenant à une classe ambige de $K_m$. Écrivons $\a_m^{\gamma-1}=(\alpha_m)$ pour un $\alpha_m\in K_m$, posons $\alpha_n= N_{K_{m\!}/\!K_n}(\alpha_m)$ et considérons l'unité $\varepsilon=N_{K_m/K}(\alpha_m) \in E_K \cap N_{K_m/K}(K_m^\times)$.
Par construction $\varepsilon$ est une puissance $\ell^m$-ième locale aux places au-dessus de $\ell$, donc une puissance $\ell^n$-ième globale, disons $\varepsilon=\eta^{\ell^n}$ pour un $\eta \in E_K$. Il suit: $N_{K_n/K}(\alpha_n/\eta)=1$, donc $\alpha_n=\eta\beta_n^{\gamma-1}$ pour un $\beta_n$ de $K_n^\times$. Et l'idéal $\a_n/\beta_n$ est un idéal ambige de $K_n$, produit comme tel d'un idéal construit sur les premiers au-dessus de $\ell$ et d'un idéal étendu de $K$ donc principal dans $K_n$.
D'où $(xii)$.
\end{itemize}\medskip

Sous les hypothèses du Théorème \ref{CrGreenberg}, la conjecture de Greenberg affirme en conséquence que le groupe sauvage $\,\C^{[\ell]}_K=\varprojlim \,\Cl^{[\ell]}_{K_n}$ ne peut être trivial sans qu'on ait simultanément $\,\Cl_{K_n}^{\,\Gamma}=1$, pour tout $n\gg 0$; et donc, d'après le Corollaire \ref{CGenres}: $|\,\T_K^{bp}|=|{}^\Gamma \Cl_{K_n}|=1$. En particulier:

\begin{Sco}\label{rationnel}
Sous les hypothèses du Théorème \ref{CrGreenberg}, la condition $\,\C_K^{[\ell]}=1$ ne peut être satisfaite sous la conjecture de Greenberg que si le module de Bertrandias-Payan $\,\T_K^{bp}$ du corps $K$ est trivial.
\end{Sco}

Si $\ell$ est impair, il suit de là, comme observé par \cite{Gra5}, qu'un tel corps est $\ell$-rationnel (cf. aussi \cite{J56}; ainsi que \cite{GJ,J56,JN, Mov, MN} sur les relations entre $\ell$-rationalité et $\ell$-régularité).

\newpage
%%%%%%%%%%%%%%%%%%%%%%%%%%%%%%%%%%%%%%%%%%%%%%%%%%%%%%%%%%%
\subsection{Schéma abélien des principales pro-$\ell$-extensions $\ell$-ramifiées}
%%%%%%%%%%%%%%%%%%%%%%%%%%%%%%%%%%%%%%%%%%%%%%%%%%%%%%%%%%%

En vue de compléter le schéma galoisien donné dans la section \ref{2.1}, introduisons le compositum $K_\infty K^{nr}$ de la $\Zl$-extension cyclotomique $K_\infty$ et de la $\ell$-extension abélienne non-ramifiée $\infty$-décomposée $K^{nr}$ de $K$. Son groupe de normes est ainsi:\smallskip

\centerline{$\wJ_K \cap \,\U_{K_\ell}\prod_{\p\nmid\ell}\mu_\p\,\R_K=\,\U^*_{K_\ell}\prod_{\p\nmid\ell}\mu_\p\,\R_K$.}\smallskip

\noindent Comparons $K_\infty K^{nr}$ à la pro-$\ell$-extension localement cyclotomique $K^{lc}$ de groupe de normes:\smallskip

\centerline{$\wU_K\R_K=\,\wU_{K_\ell}\prod_{\p\nmid\ell}\mu_\p\,\R_K$.}\smallskip

\noindent De l'égalité:\qquad
$(\,\wU^{\phantom{*}}_{K_\ell}\prod\mu_\p^{\phantom{*}}\,\R^{\phantom{*}}_K)(\,\U^*_{K_\ell} (\prod\mu_\p^{\phantom{*}}\, \R^{\phantom{*}}_K)
=\,\wU^{\phantom{*}}_{K_\ell}\,\U^*_{K_\ell}\prod\mu_\p^{\phantom{*}}\,\R^{\phantom{*}}_K=\R^*_{K_\ell}\prod\mu_\p^{\phantom{*}}\,\R^{\phantom{*}}_K$,\smallskip

\noindent nous concluons que l'intersection $K^{lc}\cap K_\infty K^{nr}$ est précisément le compositum $K_\infty K^{\ell d}$, où $K^{\ell d}$ est la $\ell$-extension abélienne non-tamifiée $\ell\infty$-décomposée maximale de $K$.\smallskip
 
Pour déterminer le compositum $K^{lc}K_\infty^{nr}$, formons l'intersection $(\,\wU^{\phantom{*}}_{K_\ell}\prod\mu_\p^{\phantom{*}}\,\R^{\phantom{*}}_K) \cap (\,\U^*_{K_\ell} (\prod\mu_\p^{\phantom{*}}\, \R^{\phantom{*}}_K)$.

\begin{Lem}
Notons $\wE'_{K_\ell}$ l'image canonique du groupe des $\ell$-unités $\E'_K$ dans le groupe des unités logarithmiques semi-locales $\,\wU_{K_\ell}$, c'est-à-dire la pré-image dans $\,\wU_{K_\ell}$ du sous-groupe $s^{\phantom{*}}_\ell(\E'_K)\,\U^*_{K_\ell}/\U^*_{K_\ell}$ de $\R^*_{K_\ell}/\U^*_{K_\ell}\simeq\,\wU^{\phantom{*}}_{K_\ell}/\mu^{\phantom{*}}_{K_\ell}$.
 Écrivons de même $\E'_{K_\ell}{\!\!\!\!}^*$ l'image de $\E'_K$ dans $\,\U_{K_\ell}^*$, c'est-à-dire la pré-image dans $\,\U_{K_\ell}^*$ du sous-groupe $s^{\phantom{*}}_\ell(\E'_K)\,\wU^{\phantom{*}}_{K_\ell}/\,\wU^{\phantom{*}}_{K_\ell}$ de $\R^*_{K_\ell} /\,\wU^{\phantom{*}}_{K_\ell} \simeq \,\U^*_{K_\ell} /\mu^{\phantom{*}}_{K_\ell}$.
Il vient alors:\smallskip

\centerline{$(\,\wU^{\phantom{*}}_{K_\ell}\prod\mu_\p^{\phantom{*}}\,\R^{\phantom{*}}_K) \cap (\,\U^*_{K_\ell} (\prod\mu_\p^{\phantom{*}}\, \R^{\phantom{*}}_K)=\,\wE'_{K_\ell}\prod\mu_\p^{\phantom{*}}\R^{\phantom{*}}_K=\,\E'_{K_\ell}{\!\!\!\!}^*\prod\mu_\p^{\phantom{*}}\R^{\phantom{*}}_K$.}\smallskip

\noindent Et le groupe de Galois $\Gal(K^{bp}/K^{lc}K^{nr})$ s'identifie au quotient $ \V^{\phantom{*}}_K \underset{^\text{déf}}{=}\, \E'_K/\wE^{\phantom{*}}_K\E^{\phantom{*}}_K$.
\end{Lem}

\Preuve Le calcul de l'intersection est immédiat. Il vient, par exemple:\smallskip

\centerline{$(\,\wU^{\phantom{*}}_{K_\ell}\prod\mu_\p^{\phantom{*}}\,\R^{\phantom{*}}_K) \cap (\,\U^*_{K_\ell} (\prod\mu_\p^{\phantom{*}}\, \R^{\phantom{*}}_K)=\big(\,\wU^{\phantom{*}}_{K_\ell} \cap \prod\mu_\p^{\phantom{*}}\,\E'_K\big)\prod\mu_\p^{\phantom{*}}\,\R^{\phantom{*}}_K.$}\smallskip

Intéressons-nous donc au groupe $G=\Gal(K^{bp}/K^{lc}K^{nr})\simeq \,\wE'_{K_\ell}\prod\mu_\p^{\phantom{*}}\R^{\phantom{*}}_K \,/\,\prod\mu_\p^{\phantom{*}}\,\R^{\phantom{*}}_K$, quotient des groupes de normes respectivement attachés à $K^{lc}K^{nr}$ et à $K^{bp}$. Il vient, comme annoncé:\smallskip

\centerline{$G \simeq \wE'_{K_\ell} \,/\,(\wE'_{K_\ell}\cap\prod\mu_\p^{\phantom{*}}\,\R^{\phantom{*}}_K) \simeq  \wE'_{K_\ell} / \wE^{\phantom{*}}_{K_\ell} \simeq  \E'_{K_\ell}\U^*_{K_\ell}/ \wE^{\phantom{*}}_{K_\ell}\U^*_{K_\ell} \simeq  \E'_{K_\ell} /  \wE^{\phantom{*}}_{K_\ell} \E^{\phantom{*}}_{K_\ell} \simeq \E'_K /  \wE^{\phantom{*}}_K \E^{\phantom{*}}_K=\V^{\phantom{*}}_K$.}
\bigskip

En résumé, les extensions considérées prennent ainsi place dans le diagramme suivant où sont figurés les corps et les groupes de Galois:

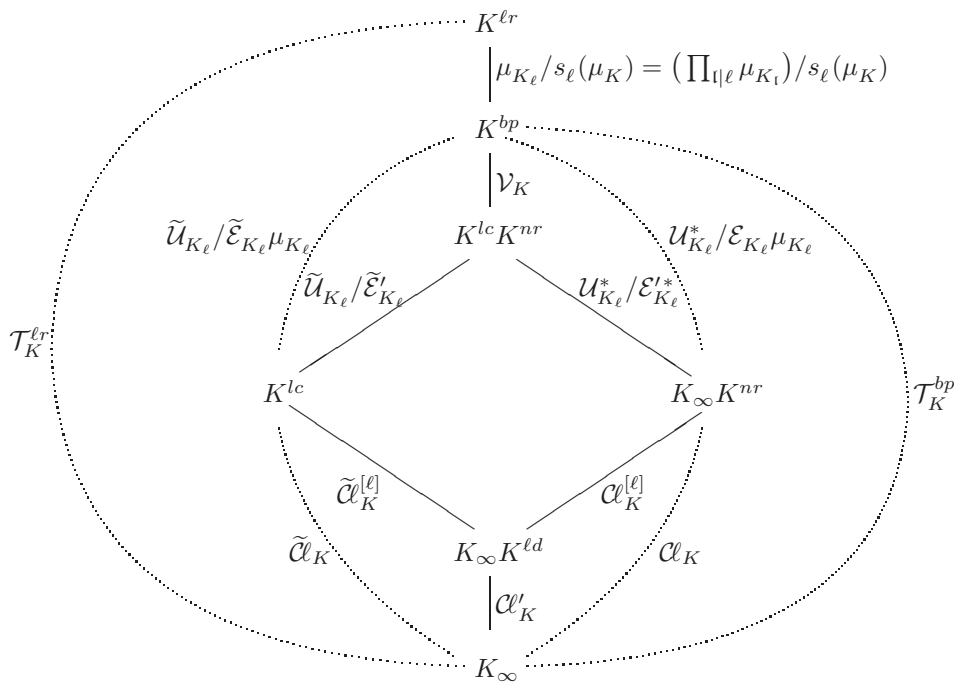
\begin{center}
\setlength{\unitlength}{2pt}
\begin{picture}(104,140)(0,0)

\thinlines 
\put(50,130){$K^{\ell r}$}
\put(50,110){$K^{bp}$}
\put(46,90){$K^{lc}K^{nr}$}
\put(10,60){$K^{lc}$}
\put(87,60){$K_\infty K^{nr}$}
\put(46,30){$K_\infty K^{\ell d}$}
\put(50,8){$K_\infty$}

\bezier{60}(46,110)(18,100)(13,70)
\bezier{60}(56,110)(87,100)(93,70)
\bezier{60}(13,55)(18,30)(46,12)
\bezier{60}(93,55)(88,30)(60,12)

\bezier{100}(60,112)(130,108)(132,60)
\bezier{100}(60,10)(130,12)(132,60)
\put(133,60){$\T^{bp}_K$}

\bezier{100}(48,132)(-28,130)(-30,70)
\bezier{100}(48,10)(-28,12)(-30,70)
\put(-38,70){$\T^{\ell r}_K$}

\put(53,127){\line(0,-1){10}}
\put(53,107){\line(0,-1){10}}
\put(53,27){\line(0,-1){10}}
\put(49,87){\line(-3,-2){32}}
\put(58,87){\line(3,-2){33}}
\put(50,36){\line(-3,2){35}}
\put(60,36){\line(3,2){33}}

\put(54,122){$\mu_{K_\ell}^{\phantom{*}}/s^{\phantom{*}}_\ell(\mu_K^{\phantom{*}})=\big(\prod_{\l|\ell}\mu_{K_\l}^{\phantom{*}}\big)/s^{\phantom{*}}_\ell(\mu_K^{\phantom{*}})$}
\put(54,100){$\V_K$}
\put(74,41){$\Cl_K^{[\ell]}$}
\put(24,41){$\wCl_K^{[\ell]}$}
\put(-8,90){$\wU^{\phantom{*}}_{K_\ell}/\wE^{\phantom{*}}_{K_\ell}\mu_{K_\ell}^{\phantom{*}}$}
\put(87,90){$\U^*_{K_\ell}/\E_{K_\ell}^{\phantom{*}}\mu_{K_\ell}^{\phantom{*}}$}
\put(15,30){$\wCl_K$}
\put(85,30){$\Cl_K$}
\put(54,20){$\Cl'_K$}
\put(18,80){$\wU^{\phantom{*}}_{K_\ell}/\wE'_{K_\ell}$}
\put(70,80){$\U^*_{K_\ell}/\E'_{K_\ell}{\!\!\!\!}^*$}
\end{picture}

\end{center}

\newpage
%%%%%%%%%%%%%%%%%%%%%%%%%%%%%%%%%%%%%%%%%%%%%%%%%%%%%%%%%%%
\subsection{Critère de Gras dans le cas $\ell$-décomposé}
%%%%%%%%%%%%%%%%%%%%%%%%%%%%%%%%%%%%%%%%%%%%%%%%%%%%%%%%%%%

Le critère de Gras consiste à remplacer la condition $(x)$ du Théorème \ref{CGreenberg}, qui exprime la trivialité de la limite inductive des $\ell$-groupes de $\ell$-classes ambiges $\,\Cl'{\!}_{K_\infty}^{\;\Gamma}=1$, par la condition plus exigeante: $\,\Cl'{\!}_{K_n}^{\;\Gamma}=1$ pour $n \gg 0$ (cf. \cite{Gra5}, Th. 2.4). Donnons en tout de suite une traduction arithmétique:

\begin{Lem}
Lorsque les places au-dessus de $\ell$ se décomposent complètement dans $K/\QQ$, on a l'équivalence (où $E'_K$ est le groupe des $\ell$-unités de $K$):

\centerline{$\displaystyle
\,\Cl'{\!}_{K_n}^{\;\Gamma}=1\quad\Leftrightarrow\quad
\begin{cases}
\quad \Cl'_K =1,\smallskip\\
\quad \big(E'_K:E'_K\cap N_{K_{n\!}/\!K}(K_n^\times)\big)= \ell^{n(d-1)}.
\end{cases}$
.}

\end{Lem}

\Preuve Cela résulte de la formule des classes ambiges qui s'écrit ici, pour les $\ell$-classes:\smallskip

\centerline{$|\Cl'{\!}_{K_n}^{\;\Gamma}|=|\Cl'_K|\frac{\prod_{\l|\ell}e_\l(K_n/K)}{[K_n:K] \big(E'_K:E'_K\cap N_{K_{n\!}/\!K}(K_n^\times)\big)}
=|\Cl'_K|\frac{\ell^{n(d-1)}}{\big(E'_K:E'_K\cap N_{K_{n\!}/\!K}(K_n^\times)\big)}$.}\medskip

Comme observé par Gras (cf. \cite{Gra5}, Th. 3.4 \& Prop. 4.6), le critère asymptotique introduit peut se lire à n'importe quel étage de la tour $K_\infty/K$. En fait, il se lit directement dans $K$:

\begin{Th}[Critère suffisant de Gras]\label{CSGras}
Soit $K$ totalement réel qui vérifie la conjecture de Leopoldt et où les places au-dessus de $\ell$ se décomposent complètement. Les quatre conditions suivantes sont alors équivalentes:\smallskip

\begin{itemize}
\item[(a)] Le critère asymptotique suffisant est satisfait: $\,\Cl'{\!}_{K_n}^{\;\Gamma}=1$, pour tout $n \gg 0$.\smallskip

\item[(b)] Le premier étage de la tour $K_\infty/K$ vérifie le critère de Gras: $\,\Cl'{\!}_{K_1}^{\;\Gamma}=1$\smallskip

\item[(c)] On a simultanément les identités: $\,\Cl'_K=1$ et $\,\U^*_{K_\ell}/\E'_{K_\ell}{\!\!\!}^*=1$.\smallskip

\item[(d)] Le $\ell$-groupe des classes logarithmiques de $K$ est trivial: $\,\wCl_K=1$.

\end{itemize}
\end{Th}

\Preuve D'après le Lemme, nous pouvons remplacer la condition $(a)$ par la réunion des deux conditions $\,\Cl'_K=1$ et $ \big(E'_K:E'_K\cap N_{K_{n\!}/\!K}(K_n^\times)\big)= \ell^{n(d-1)}$ pour $n \gg 0$. Intéressons-nous à la seconde, que nous pouvons réécrire, après $\ell$-adification et pour un $n>0$ fixé, sous la forme:\smallskip

\centerline{$\big(\E'_K:\E'_K\cap N_{K_{n\!}/\!K}(\R^{\phantom{*}}_{K_n})\big)= \ell^{n(d-1)}$.}\smallskip

\noindent  Observons d'abord que l'intersection $\E'_K\cap N_{K_{n\!}/\!K}(\R^{\phantom{*}}_{K_n})$ contient évidemment le sous-groupe $\wE^{\phantom{*}}_K$ des normes cyclotomiques ainsi que le sous-groupe des puissances $\ell^n$-ièmes dans $\E'_K$. Par suite, puisque le quotient $\E'_K/\wE^{\phantom{*}}_K$ est un $\Zl$-module libre de dimension $d-1$ (le corps $K$ vérifiant la conjecture de Gross-Kuz'min  en vertu de la Proposition \ref{Conj}), cette égalité affirme tout simplement que l'on a:

\centerline{$\E'_K\cap N_{K_{n\!}/\!K}(\R^{\phantom{*}}_{K_n})=\wE^{\phantom{*}}_K\E'_K{\!}^{\ell^n}$.}\smallskip

\noindent Or, d'après le principe de Hasse, le groupe de gauche n'est autre que la préimage dans $\E'_K$ (par l'application de semi-localisation $s^{\phantom{*}}_\ell$) du sous-groupe normique $\U^{*\ell^n}_{K_\ell}\wU^{\phantom{*}}_{K_\ell}$ de $\R^*_{K_\ell}$ associé à $K_n$.

Ainsi, dans le plongement canonique de $\E'_K/\wE^{\phantom{*}}_K$ dans $\,\U^*_{K_\ell}\simeq\R^*_{K_\ell}/\wU^{\phantom{*}}_{K_\ell}$, les éléments de l'image $\E'_{K_\ell}{\!\!\!\!}^*$ qui sont des puissances $\ell^n$-ièmes dans $\,\U^*_{K_\ell}$ sont les puissances $\ell^n$-ièmes des éléments de $\E'_{K_\ell}{\!\!\!}^*$.\smallskip

Il vient donc: $\E'_{K_\ell}{\!\!\!}^*=\U^*_{K_\ell}$; i.e. comme annoncé: $\,\U^*_{K_\ell}/\E'_{K_\ell}{\!\!\!}^*=1$. Et la réciproque est immédiate.\smallskip

Observons au passage que l'équivalence obtenue est indépendante de la valeur de $n$: elle vaut indifféremment pour tout $n \gg 0$ ou un $n>0$ arbitraire donné et donc, en particulier, pour $n=1$.\smallskip

 Enfin, conformément au diagramme des extensions dressé dans la section précédente, la conjonction des deux égalités  $\,\Cl'_K=1$ et $\,\U^*_{K_\ell}/\E'_{K_\ell}{\!\!\!}^*=1$ traduit simplement la trivialité du groupe $\,\wCl_K$.\medskip

\Exemple Pour $k=\QQ[\sqrt{m}]$ quadratique réel,  $\ell=3$ et $m \equiv 1\;[\text{mod }3]$, sur les quelque 2256 corps obtenus pour $m < 10.000$, seuls 237 (soit environ 10\%) ont un 3-groupe des classes logarithmiques non trivial\footnote{Ce sont précisément ceux apparaissant dans les tables de Fukuda et Taya \cite{FuTa}}. Cette proportion est de 2.801/22.793 (soit environ 12 \%) pour $ m<100.000$. Elle est de 30.747/227.953 (soit environ 13,5 \%) pour $m<1.000.000$.

\newpage
%%%%%%%%%%%%%%%%%%%%%%%%%%%%%%%%%%%%%%%%%%%%%%%%%%%%%%%%%%%%%%%%%%%%%%
%%%%%%%%%%%%%%%%%%%%%%%%%%%%%%%%%%%%%%%%%%%%%%%%%%%%%%%%%%%%%%%%%%%%%%
\section{Trivialité du sous-module fini et conjecture faible}
%%%%%%%%%%%%%%%%%%%%%%%%%%%%%%%%%%%%%%%%%%%%%%%%%%%%%%%%%%%%%%%%%%%%%%
%%%%%%%%%%%%%%%%%%%%%%%%%%%%%%%%%%%%%%%%%%%%%%%%%%%%%%%%%%%%%%%%%%%%%%
%$\overset{_\leftarrow}{\E}_F=\varprojlim \E_{F_n}$

%$\overset{_\leftarrow}{\E}{}^\circ_F=\varprojlim \E^\circ_{F_n}$

Intéressons-nous maintenant aux plus grands sous-modules finis des $\Lambda$-modules considérés.\medskip

Il est bien connu que ces sous-modules s'interprètent comme limites projectives des sous-groupes de capitulation. Rappelons brièvement pourquoi. 
Dans tous les cas (i.e. pour $X_n=\,\Cl_{K_n}$ ou $\,\Cl'_{K_n}$ ou encore $\,\wCl_{K_n}$) le schéma général est le même : le $\Lambda$-module $X$ est la limite projective des $X_n$ (pour le morphismes normiques) et les $X_n$ s'obtiennent en retour pour $n\gg 0$, disons $n\ge n_\circ$, à partir de $X$ comme quotients

\centerline{$X_n=X/\frac{\omega_{n_{\phantom{\!l}}}}{\omega_{n_{\scale{.6}{\scriptstyle \circ}}}}\,Y$,}\smallskip

\noindent où $Y$ est un sous-module convenable de $X$ qui contient $\omega_{n_\circ}X$ comme sous-module de cotype fini. Dans le cas des $\ell$-groupes de classes, $Y/\omega_{n_\circ}X$ est engendré par les images des sous-groupes d'inertie attachés aux places au-dessus de $\ell$; dans le cas des $\ell$-groupes de $\ell$-classes, il est engendré par les images des sous-groupes de décomposition de ces mêmes places; et dans le cas des $\ell$-groupes de classes logarithmiques c'est encore plus simple, car on a tout simplement: $n_\circ=0$ et $Y=\omega_{n_\circ}X$.\smallskip

Dans tous les cas les morphismes d'extension $j_{n,m}:\;X_n \rightarrow X_m$, pour $m \ge n$, sont donnés par la multiplication par $\omega_m/\omega_n$. Et il vient donc:\smallskip

\centerline{$\Ker\,j_{n,m} = \{ x+\frac{\omega_{n_{\phantom{\!l}}}}{\omega_{n_{\scale{.6}{\scriptstyle \circ}}}}\,Y \in X/\frac{\omega_{n_{\phantom{\!l}}}}{\omega_{n_{\scale{.6}{\scriptstyle \circ}}}}\,Y \, | \, \frac{\omega_{m_{\phantom{\!l}}}}{\omega_n}\,x  \in \frac{\omega_{m_{\phantom{\!l}}}}{\omega_{n_{\scale{.6}{\scriptstyle \circ}}}}\,Y\}$.}\smallskip 

\noindent Soit, comme précédemment, $F$ le plus grand sous-module fini de $X$. Écrivant alors $\frac{\omega_{m_{\phantom{\!l}}}}{\omega_n}\,x=\frac{\omega_{m_{\phantom{\!l}}}}{\omega_{n_{\scale{.6}{\scriptstyle \circ}}}}\,y$ pour un $y \in Y$, on obtient  $\frac{\omega_{m_{\phantom{\!l}}}}{\omega_n}\,(x-\frac{\omega_{n_{\phantom{\!l}}}}{\omega_{n_{\scale{.6}{\scriptstyle \circ}}}} y)=0$, i.e. $x-\frac{\omega_{n_{\phantom{\!l}}}}{\omega_{n_{\scale{.6}{\scriptstyle \circ}}}} y \in F$, pour $m \gg n\gg 0$, puisque le polynôme $\frac{\omega_{m_{\phantom{\!l}}}}{\omega_n}$ est étranger au polynôme caractéristique du $\Lambda$-module $X$ pour $n \gg 0$ et que par ailleurs il annule bien $F$ pour $m \gg n$. En résumé, il vient:\smallskip

\centerline{$\Ker\,j_{n,m} \underset{m \gg n \gg 0}{=} (F + \frac{\omega_{n_{\phantom{\!l}}}}{\omega_{n_{\scale{.6}{\scriptstyle \circ}}}}\,Y ) / \frac{\omega_{n_{\phantom{\!l}}}}{\omega_{n_{\scale{.6}{\scriptstyle \circ}}}}\,Y 
\simeq F/(F\cap  \frac{\omega_{n_{\phantom{\!l}}}}{\omega_{n_{\scale{.6}{\scriptstyle \circ}}}}\,Y)
\underset{n \gg n_\circ}{=} F $;}\smallskip

\noindent de sorte que le sous-module fini $F$ de $X$ s'identifie ultimement au sous-groupe de capitulation:\smallskip

\centerline{$\Ker\,j_{n,\infty}:\;X_n \rightarrow X_\infty=\varinjlim X_m$.}

\begin{Th}
Soient $K$ un corps de nombres arbitraire, $\ell$ un nombre premier.  et $K_\infty=\cup_{n\in\NN}K_n$ la $\Zl$-extension cyclotomique de $K$. On a alors inconditionnellement les interprétations suivantes:
\begin{itemize}
\item[(i)] Le sous-module fini $\,\F_K$ de la limite projective des $\ell$-groupes de classes $\,\C_K=\varprojlim \,\Cl_{K_n}$ s'identifie ultimement aux noyaux des morphismes d'extension $\,Cap_n=\Ker (\,\Cl_{K_n}\rightarrow\,\Cl_{K_\infty})$.
\item[(ii)] Le sous-module fini $\,\F'_K$ de la limite projective des $\ell$-groupes de $\ell$-classes $\,\C'_K=\varprojlim \,\Cl'_{K_n}$ s'identifie ultimement aux noyaux des morphismes d'extension $\,Cap'_n=\Ker (\,\Cl'_{K_n}\rightarrow\,\Cl'_{K_\infty})$.
\item[(iii)] Le même sous-module $\,\F'_K$ s'identifie encore ultimement aux noyaux des morphismes d'extension entre groupes de classes logarithmiques $\,\wCap_n=\Ker (\,\wCl_{K_n}\rightarrow\,\wCl_{K_\infty})$.
\end{itemize}
\end{Th}

Il suit de là que le sous-module fini $F$ de $X$ est trivial si et seulement si les morphismes d'extension $j_{n,\infty}:\;X_n \rightarrow X_\infty$ sont ultimement injectifs. Mais il est possible d'être plus précis:

\begin{Sco}
Soit $K$ est un corps totalement réel qui vérifie la conjecture de Leopoldt en $\ell$. Alors:\smallskip
\begin{itemize}
\item[(i)] Si les places au-dessus de $\ell$ sont totalement ramifiées dans $K_\infty/K$, on a l'équivalence:\

\centerline{$\F_K=1 \quad \Leftrightarrow \quad Cap_n=1$, pour chaque $n\in\NN$.}\smallskip

\item[(ii)] Si aucune des places au-dessus de $\ell$ ne se décompose dans $K_\infty/K$, on a l'équivalence:\

\centerline{$\F'_K=1 \quad \Leftrightarrow \quad Cap'_n=1$, pour chaque $n\in\NN$.}\smallskip

\item[(iii)] Et inconditionnellement, on a l'équivalence:\

\centerline{$\F'_K=1 \quad \Leftrightarrow \quad \wCap_n=1$, pour chaque $n\in\NN$.}
\end{itemize}
\end{Sco}

\Preuve La conjecture de Leopoldt assure que le polynôme caractéristique du $\Lambda$-module $X$ est étranger aux $\omega_n$ (cf. Proposition \ref{Conj}).
Et, dans les deux premiers cas, l'hypothèse faite permet de prendre $n_\circ=0$. Cela étant, dans les trois cas, si $F$ est nul, la condition de capitulation obtenue  plus haut $\frac{\omega_{m_{\phantom{\!l}}}}{\omega_n}\,(x-\frac{\omega_n}{\omega_{n_{\scale{.6}{\scriptstyle \circ}}}} y)=0$ donne directement $x=\frac{\omega_{n_{\phantom{\!l}}}}{\omega_{n_{\scale{.6}{\scriptstyle \circ}}}} y \in \frac{\omega_{n_{\phantom{\!l}}}}{\omega_{n_{\scale{.6}{\scriptstyle \circ}}}} Y$, i.e. $\Ker\,j_{n,m} =0$, pour $m \ge n\ge 0$.\smallskip

Et les morphismes d'extension $j_{n,\infty}$ sont bien injectifs pour tout $n\in\NN$.

\newpage
%%%%%%%%%%%%%%%%%%%%%%%%%%%%%%%%%%%%%%%%%%%%%%%%%%%%%%%%%%%
\subsection{Sous-modules fini et sous-module sauvage}\medskip
%%%%%%%%%%%%%%%%%%%%%%%%%%%%%%%%%%%%%%%%%%%%%%%%%%%%%%%%%%%

Pour étudier la trivialité du sous-module $\,\F_K$, précisons en différentes formulations:

\begin{Prop}\label{CGFaible}
Soit $K$ un corps de nombres totalement réel; $K_\infty$ sa $\Zl$-extension cyclotomique; $K_{d_0^{\phantom{*}}}$ le compositum des sous-corps de décomposition des places au-dessus de $\ell$ dans $K_\infty/K$; et $K_{r_0^{\phantom{*}}}$ le compositum de leurs sous-corps d'inertie (au sens habituel). Si le corps $K_{d_0^{\phantom{*}}}$ satisfait la conjecture de Leopoldt (pour le premier $\ell$), la trivialité du sous-module fini $\,\F_K$ de $\,\C_K$ entraîne celle du sous-module sauvage $\,\C_K^{[\ell]}$, i.e. l'identité: $\,\Cl^{\phantom{*}}_{K_n}=\,\Cl'_{K_n}$, pour $n \ge r_0^{\phantom{*}}$.\smallskip

En d'autres termes, on a alors l'équivalence:\quad $\,\F_K=1 \quad\Leftrightarrow\quad \C^{[\ell]}_K=1 \quad\&\quad \F'_K=1$.
\end{Prop}

\Preuve Observons d'abord que  les morphismes normiques $ N_{K_{m\!}/\!K_n}:\,\Cl^{[\ell]}_{K_m}\rightarrow\,\Cl^{[\ell]}_{K_n}$ entre sous-groupes sauvages sont surjectifs pour $m \ge n \ge r_0^{\phantom{*}}$, puisque les places au-dessus de $\ell$ sont totalement ramifiées dans  $K_\infty/K_{r_0^{\phantom{*}}}$. Ainsi, il revient donc au même d'affirmer que les sous-groupes sauvages respectifs des corps $K_n$ sont triviaux pour $n \ge r_0^{\phantom{*}}$ ou que leur limite projective l'est.\smallskip
\begin{itemize}

\item Supposons $\,\F_K=1$. Si le corps $K_{d_0^{\phantom{*}}}$ vérifie la conjecture de Leopoldt (en $\ell$), le sous-groupe sauvage de $\,\C_K$, qui est fini d'après le Scolie \ref{Classes sauvages}, est donc trivial.
On conclut $\,\Cl^{\phantom{*}}_{K_n}=\,\Cl'_{K_n}$, pour $n \ge r_0^{\phantom{*}}$; et il suit en particulier: $\,\C'_K=\,\C_K$ puis $\,\F'_K=1=\,\C^{[\ell]}_K$, comme annoncé.\smallskip

\item Inversement, supposons  $\,\F'_K=1=\,\C_K^{[\ell]}$. La seconde condition nous donne: $\,\Cl^{[\ell]}_{K_n}=1$ pour $n \ge r_0^{\phantom{*}}$, i.e. $\,\Cl^{\phantom{*}}_{K_n}=\,\Cl'_{K_n}$; et il suit: $\,\C_K=\,\C'_K$; puis: $\,\F^{\phantom{*}}_K=\,\F'_K=1$.
\end{itemize}\medskip

Examinons maintenant les deux conditions de trivialité à droite: $ \F'_K=1$ et $\, \C^{[\ell]}_K=1$.

\begin{Prop}\label{Kuz}
Si le corps $K$ vérifie la conjecture de Gross-Kuz'min (pour le premier $\ell$), le sous-module fini $\,\F'_K$ du module de Kuz'min-Tate $\,\C'_K$ est trivial si et seulement si le $\ell$-groupe des unités logarithmiques $\,\wE_K$ coïncide avec le sous-groupe normique $\,\wE_K^\nu=\bigcap_{n\in\NN}N_{K_n/K}(\E'_{K_n})$:

\centerline{$\F'_K=1 \quad \Leftrightarrow \quad \wE_K^{\phantom{*}} = \wE_K^\nu$.}
\end{Prop}

\Preuve C'est une conséquence facile de l'isomorphisme de Kuz'min: $\,\C'_K{}^{\!\Gamma} \simeq \wE_K/ \wE_K^\nu$ (cf. e.g. \cite{J55}, Th. 17). Sous la conjecture de Gross-Kuz'min, il vient, en effet: $\,\C'_K{}^{\!\Gamma}\subset\,\F'_K$; d'où les équivalences:\smallskip

\centerline{$\,\F'_K=1  \Leftrightarrow \,\F'_K{}^{\!\Gamma}=1 \Leftrightarrow \,\C'_K{}^{\!\Gamma}=1 \Leftrightarrow \wE_K /\wE_K^\nu=1$.}

\begin{Prop}\label{Ram}
Si les places sauvages sont totalement ramifiées dans $K_\infty/K$, on a l'implication:\smallskip

\centerline{$\C_K^{[\ell]}=1 \quad \Rightarrow \quad \E'_K=\,\wE_K^\nu\,\E_K^{\phantom{*}}$.}\smallskip

\noindent Et $\,\E'_K$ est alors engendré par le sous-groupe normique $\,\wE_K^\nu$ et le groupe des unités $\,\E_K^{\phantom{*}}$.
\end{Prop}

\Preuve Si les places sauvages sont totalement ramifiées dans $K_\infty/K$, la norme arithmétique envoie le sous-module sauvage $\D_{K_n}^{[\ell]}$  du groupe des diviseurs de $K_n$ sur son homologue dans $K$. L'égalité entre diviseurs $\D_K^{[\ell]}=N_{K_n/K}(\D_{K_n}^{[\ell]})$ donne ainsi $\,\E'_F=N_{K_n/K}(\E'_{K_n})\E_K^{\phantom{*}}$, lorsque ces diviseurs sauvages sont supposés principaux. Faisant varier $n$, nous obtenons alors, par compacité de $\,\E_F$:\smallskip

\centerline{$\E'_K\,=\,\bigcap_{n\in\NN}\big(N_{K_n/K}(\E'_{K_n})\E_K^{\phantom{*}}\big)\,=\,\big(\bigcap_{n\in\NN}N_{K_n/K}(\E'_{K_n})\big)\E_K^{\phantom{*}}\,=\,\wE_K^\nu\,\E_K^{\phantom{*}}$.}\smallskip

\begin{Th}\label{CrSauvage}
Soient $K$ un corps totalement réel, $d=[K:\QQ ]$ son degré et $\ell$ un premier complètement décomposé dans $K$. Sous la conjecture de Leopodt pour $\ell$ dans $K$, on a l'équivalence:

\centerline{$\,\F^{\phantom{*}}_K=1 \quad\Leftrightarrow\quad \C^{[\ell]}_K=1$.}
\end{Th}

\noindent{\em Preuve} (voir aussi \cite{Ng2}, Lem. 2.6). Supposons $\,\C^{[\ell]}_K=1$. De l'égalité $\,\E'_K=\,\wE_K^\nu\,\E_K$, avec ici $\;\dim_{\Zl}\wE_K^\nu=\dim_{\Zl}\wE_K^{\phantom{*}}=d$ (puisque $K$, qui vérifie la conjecture de Leopoldt, vérifie aussi celle de Gross-Kuz'min) et $\dim_{\Zl}\E_K^{\phantom{*}}=d-1$, donc $\;\dim_{\Zl}\wE_K^\nu+\dim_{\Zl}\E_K^{\phantom{*}}=\dim_{\Zl}\E'_K$, nous concluons que la décomposition $\,\E'_K=\,\wE_K^\nu\,\E_K$ est directe aux racines de l'unité près.\par
En particulier $\,\E'_K/\,\wE_K^\nu\,\simeq\,\E_K^{\phantom{*}}/(\,\E_K^{\phantom{*}}\cap\,\wE_K^\nu)\,\simeq \,\E_K^{\phantom{*}}/\mu_K^{\phantom{*}}$ est $\Zl$-libre; et il suit: $\,\wE_K^{\phantom{*}} = \wE_K^\nu$, c'est-à-dire $\,\F'_K=1$; d'où finalement: $\,\F^{\phantom{*}}_K=1$, comme annoncé.

\newpage

%%%%%%%%%%%%%%%%%%%%%%%%%%%%%%%%%%%%%%%%%%%%%%%%%%%%%%%%%%%
\subsection{Étude du cas abélien semi-simple et $\ell$-décomposé: préliminaires}
%%%%%%%%%%%%%%%%%%%%%%%%%%%%%%%%%%%%%%%%%%%%%%%%%%%%%%%%%%%

Nous supposons désormais que $K$ est un corps abélien réel de groupe de Galois $\Delta=\Gal(K/\QQ)$ et que $\ell$ est un nombre premier {\em impair} complètement décomposé dans $K$ et ne divisant pas $d=[K:\QQ]$. Rappelons que l'hypothèse de complète décomposition implique en particulier que les places $\l$ de $K$ au-dessus de $\ell$ se ramifient totalement dans la $\Zl$-tour $K_\infty/K$.\smallskip
\smallskip

Dans ce contexte, l'algèbre $\Zl[\Delta]$ est une algèbre semi-locale qui s'écrit comme produit direct\smallskip

\centerline{$\Zl[\Delta]=\oplus_\varphi \Zl[\Delta] e_\varphi=\oplus_\varphi Z_\varphi$}\smallskip

\noindent d'extension non-ramifiées $Z_\varphi$ de $\Zl$ indexées par les caractères $\ell$-adiques irréductibles de $\Delta$. Et les idempotents primitifs associés $e_\varphi=\frac{1}{d}\sum_{\sigma\in\Delta}\varphi(\tau)\tau^{-1}$ permettent de décomposer chaque $\Zl[\Delta]$-module comme somme directe de ses $\varphi$-composantes, lesquelles sont des $Z_\varphi$-modules.\smallskip

Cette décomposition s'étend naturellement à l'algèbre d'Iwasawa $\Lambda=\Zl[[\gamma-1]]$ de sorte que\smallskip

\centerline{$\Lambda[\Delta]=\oplus_\varphi \Lambda[\Delta] e_\varphi=\oplus_\varphi \Lambda_\varphi$}\smallskip

\noindent est le produit direct des algèbres de séries formelles $\Lambda_\varphi=Z_\varphi[[\gamma-1]]$ construites sur les $Z_\varphi$.
\medskip

Introduisons maintenant les $\ell$-adifiés $\,\E^\circ_{K_n}$ des groupes d'unités circulaires à la Sinnott (cf. \cite{Si}):
Notons  $\,\overset{\!\!_{\leftarrow}}{\E'_K}=\varprojlim \E'_{K_n}$, $\,\overset{_\leftarrow}{\E}\!_K=\varprojlim \E_{K_n}$ et $\,\overset{_\leftarrow}{\E}{}_{\!K}^\circ=\varprojlim \E^\circ_{K_n}$ les limites projectives respectives pour les applications normes des groupes de $\ell$-unités, d'unités et d'unités circulaires dans la $\Zl$-extension cyclotomique $K_\infty/K$.
Écrivons enfin $\,\wE_{K_n}^\nu=\bigcap_{m\ge n} N_{K_{m\!}/\!K_n}(\E'_{K_m})=\bigcap_{m\ge n} N_{K_{m\!}/\!K_n}(\wE_{K_m})$ puis $\,\E_{K_n}^\nu=\bigcap_{m\ge n} N_{K_{m\!}/\!K_n}(\E_{K_m})$ et $\,\E_{K_n}^{\circ\, \nu}=\bigcap_{m\ge n} N_{K_{m\!}/\!K_n}(\E^{\circ}_{K_m})$ les sous-groupes normiques associés.

Il est bien connu que  $\,\overset{_{\!\leftarrow}}{\E'_K}$ comme  $\,\overset{_\leftarrow}{\E}_K$ sont $\Lambda$-libres, en vertu d'un théorème de Kuz'min (\cite{Kuz}, Th. 7.2) généralisé par Greither (\cite{Grt}, Th. 1). Et, d'après Belliard (\cite{Be}), le même résultat vaut pour les unités circulaires dès lors que celles-ci vérifient la descente galoisienne asymptotique dans la tour $K_\infty/K$, ce qui a bien lieu dans le cas semi-simple considéré ici. Il suit de là et du Lemme 3.3 de \cite{BB} que chacune des trois limites projectives ci-dessus est un $\Lambda[\Delta]$-module libre de dimension 1.

Faisons choix d'un $\Lambda[\Delta]$-base $\eta=(\eta_n)_n$ de  $\,\overset{\!\!_{\leftarrow}}{\E'_K}$. Observons que les $\eta_n$ forment un système cohérent de $\ell$-unités et qu'on a en particulier: $\eta_n\in \,\wE_{K_n}^\nu \subset \,\wE_{K_n}^{\phantom{*}}$; en d'autres termes, que les $\eta_n$ sont des unités logarithmiques, normes d'unités logarithmiques à chaque étage de la tour $K_\infty/K_n$.

La suite exacte courte canonique $1 \rightarrow \E_{K_n}^{\phantom{*}}\rightarrow \E'_{K_n} \rightarrow \P_{K_n}^{[\ell]} \rightarrow 1$, où $ \P_{K_n}^{[\ell]}$ désigne le sous-groupe principal du module des diviseurs de $K_n$ construits sur les places au-dessus de $\ell$, donne à la limite:
\begin{equation}\label{Limproj}
1 \rightarrow \overset{_\leftarrow}{\E}\!_K^{\phantom{*}}=\varprojlim \E_{K_n}^{\phantom{*}} \rightarrow \overset{_\leftarrow}{\E}{}'_K=\varprojlim \E'_{K_n}\simeq \Lambda[\Delta] \rightarrow  \overset{_\leftarrow}{\P}\!_K^{\,[\ell]}=\varprojlim\P_{K_n}^{[\ell]} \rightarrow 1.\end{equation}
\noindent Or, les morphismes normiques $\D_{K_m}^{[\ell]} \rightarrow\D_{K_n}^{[\ell]}$ et $\P_{K_m}^{[\ell]} \rightarrow\P_{K_n}^{[\ell]}$ sont ultimement surjectifs, les premiers du fait de la totale ramification sauvage , les seconds parce que les quotients correspondants stationnent, comme expliqué dans la section 1.1, puisque le corps $K$, supposé abélien, vérifie la conjecture de Leopoldt. La limite projective à droite est donc isomorphe à $\Zl[\Delta]$ de sorte que le noyau à gauche s'identifie au sous-module de l'algèbre $\Lambda[\Delta]$ engendré par l'élément $\omega=\gamma-1$.

En résumé $\, \overset{_\leftarrow}{\E}\!_K^{\phantom{*}}$ est le $\Lambda[\Delta]$-module libre engendré par la famille $\varepsilon=\eta^\omega$. Ainsi:

\begin{Th}\label{Ucirc}
Soient $\ell$ un nombre premier impair, $K$ un corps abélien réel de degré $d$ étranger à $\ell$ dans lequel $\ell$ est complètement décomposé, $K_\infty = \bigcup_{n\in\NN}K_n$ la $\Zl$-extension cyclotomique de $K$ et $\Delta$ le groupe de Galois $\Gal(K/\QQ)$. Alors:
\begin{itemize}
\item[(i)] La limite projective $\,\overset{\!\!_{\leftarrow}}{\E'_K}=\varprojlim \E'_{K_n}=\varprojlim \wE_{K_n}^{\phantom{*}}$ est un $\Lambda[\Delta]$-module isomorphe à $\Lambda[\Delta]$.
\item[(ii)]  La limite projective $\,\overset{_\leftarrow}{\E}\!_K^{\phantom{*}}=\varprojlim \E_{K_n}^{\phantom{*}}=\varprojlim \E_{K_n}^\nu$ s'identifie au sous-module engendré par $\omega$.\smallskip
\item[(iii)] Il existe une famille $(\rho_\varphi)_\varphi$ de polynômes distingués $\rho_\varphi \in Z_\varphi[\omega]$, indexée par les caractères $\ell$-adiques irréductibles du groupe $\Delta$, et un élément $\rho=\sum_\varphi \rho_\varphi e_\varphi\in\Zl[\omega][\Delta]$ tels que la limite projective $\,\overset{_\leftarrow}{\E}\!_K^{\,\circ}=\varprojlim \E_{K_n}^{\circ}=\varprojlim \E_{K_n}^{\circ\,\nu}$ s'identifie au sous-module $\rho\,\omega\Lambda[\Delta]$.
\end{itemize}
\end{Th}

\Preuve La $\varphi$-composante de $\,\overset{_\leftarrow}{\E}\!_K^{\,\circ}$ est, en effet, un sous-module libre du $\Lambda_\varphi$-module $\,\overset{_\leftarrow}{\E}\!_K^{\,e_\varphi}\simeq \Lambda_\varphi$, donc de la forme $\,\overset{_\leftarrow}{\E}\!_K^{\,\rho_\varphi} \simeq \rho_\varphi\Lambda_\varphi$, pour un polynôme distingué convenable de l'algèbre $Z_\varphi[\omega] \subset \Lambda_\varphi$.

\newpage

%%%%%%%%%%%%%%%%%%%%%%%%%%%%%%%%%%%%%%%%%%%%%%%%%%%%%%%%%%%
\subsection{Groupe des classes circulaires et lien avec la conjecture de Greenberg}\medskip
%%%%%%%%%%%%%%%%%%%%%%%%%%%%%%%%%%%%%%%%%%%%%%%%%%%%%%%%%%%

Classiquement les $\ell$-groupe de classes circulaires sont définis comme suit:

\begin{Def}
Étant donné un corps abélien réel $K$ et un nombre premier impair $\ell$, on appelle $\ell$-groupe des classes circulaires du corps $K$ le quotient\smallskip

\centerline{$\Cl^{\,\circ}_K = \,\E_K^{\phantom{*}}/\,\E^\circ_K$}\smallskip

\noindent du $\ell$-adifié $\,\E_K^{\phantom{*}}$ du groupe des unités de $K$ par celui du sous-groupe des unités circulaires de Sinnott.\smallskip

Nous notons $\,\C^{\,\circ}_K = \varprojlim\,\Cl^{\,\circ}_{K_n}\simeq\,\overset{_\leftarrow}{\E}\!_K^{\phantom{*}} / \,\overset{_\leftarrow}{\E}{\!}^{\,\circ}_K$ la limite projective pour la norme des $\ell$-groupes de classes circulaires attachés aux étages finis $K_n$ de la $\Zl$-tour cyclotomique $K_\infty/K$.
\end{Def}

Il est bien connu que les $\ell$-groupes de classes circulaires $\,\Cl^{\,\circ}_{K_n}$ sont finis et que, pour $\ell\ne 2$, leurs ordres respectifs coïncident avec ceux des $\ell$-groupes des classes $\,\Cl_{K_n}$. 
Le Théorème de Mazur-Wiles affirme alors que leur limite projective  $\,\C^{\,\circ}_K$ est encore un $\Lambda$-module noethérien et de torsion, qui a même polynôme caractéristique que  $\,\C^{\phantom{*}}_K$.  De plus, sous l'hypothèse de semi-simplicité, son sous-module fini $\,\F_K^{\,\circ}$ est trivial, puisque  $\,\C^{\,\circ}_K$ est quotient de deux $\Lambda[\Delta]$-modules libres de dimension 1. La conjecture de Greenberg pour le corps $K$ et le premier $\ell$ s'écrit donc tout simplement:$\;\C^{\,\circ}_K=1$.\smallskip

Plus précisément, dans le cas semi-simple, le Théorème de Mazur-Wiles nous donne:

\begin{Prop}
Sous les hypothèses du Théorème \ref{Ucirc}, pour chaque caractère $\ell$-adique de $\Delta$ le polynôme distingué $\rho_\varphi$  est le polynôme caractéristique de la $\varphi$-composante du $\Lambda[\Delta]$-module $\,\C_K$.\par
En particulier, la conjecture de Greenberg postule que tous les $\rho_\varphi$ sont égaux à 1.
\end{Prop}

Revenons maintenant sur le suite exacte (\ref{Limproj}). Faisant agir l'élément $\omega_n=\gamma^{\ell^n}-1$ et formant la suite exacte du serpent, nous obtenons la ligne supérieure du diagramme commutatif:

\begin{displaymath}
\xymatrix{
1 \ar[r] 	&	\overset{_\leftarrow}{\P}\!_K^{\,[\ell]}  \ar[r] 
		&	\overset{_\leftarrow}{\E}\!_K^{\phantom{*}}/\overset{_\leftarrow}{\E}\!_K^{\;\omega_n}  \ar@{->>}[d]^{N_n} \ar[r] 
		&	\overset{_\leftarrow}{\E}{}'_K/\overset{_\leftarrow}{\E}{}'\!_K^{\;\omega_n} \ar@{^{(}->>}[d]^{N_n} \ar[r] 
		&  \overset{_\leftarrow}{\P}{\!}_K^{\,[\ell]} \ar[r] \ar[d]^{N_n}
		& 1 \\
		& 1 \ar[r] 	&  \E_{K_n}^{\nu} \ar[r] &  \wE_{K_n}^\nu \ar[r] & \P_{K_n}^{[\ell]} &%\ar[r] & 1
		}
\end{displaymath}
qui, par comparaison avec la suite inférieure, nous fournit la suite exacte courte:
\begin{displaymath}
1 \longrightarrow \overset{_\leftarrow}{\P}\!_K^{\,[\ell]}  \longrightarrow\overset{_\leftarrow}{\E}\!_K^{\phantom{*}}/\overset{_\leftarrow}{\E}\!_K^{\;\omega_n} \longrightarrow \E_{K_n}^\nu  \longrightarrow  1,
\end{displaymath}
puisque la flèche verticale au centre du diagramme est un isomorphisme. 
De $\overset{_\leftarrow}{\E}\!_K^{\phantom{*}}= \varepsilon^{\Lambda[\Delta]} \simeq \Lambda[\Delta]$, nous concluons par le  Théorème \ref{Ucirc} que le groupe au centre est isomorphe au quotient $\Lambda[\Delta]/\omega_n \Lambda[\Delta]$. Et le noyau à gauche, qui est isomorphe à $\Zl[\Delta]$, s'identifie au sous-quotient $\frac{\omega_n}{\omega}\Lambda[\Delta]/\omega_n \Lambda[\Delta]$ annulé par $\omega$ (puisqu'ils ont même rang et que le conoyau $\E_{K_n}^\nu$ est sans torsion). Il vient ainsi:\smallskip

\centerline{$\E_{K_n}^\nu=\varepsilon_n^{\Lambda[\Delta]}\simeq \Lambda[\Delta]/\frac{\omega_n}{\omega}\Lambda[\Delta]$,}\smallskip

\noindent avec $\varepsilon_n=N_n(\varepsilon)=N_n(\eta^{\,\omega})=\eta_n^{\,\omega}$. Il suit: $\,\E_{K_n}^\nu=\wE_{K_n}^{\nu\,\omega}$; et, en particulier: $\,\E_K^\nu=1$.

Enfin, par restriction au sous-module des unités circulaires $\,\overset{_\leftarrow}{\E}{\!}^{\,\circ}_K= \varepsilon^{\rho\Lambda[\Delta]}$, nous obtenons:\smallskip

\centerline{$\E_{K_n}^{\circ\nu}=\varepsilon_n^{\rho\Lambda[\Delta]}
\simeq \rho\Lambda[\Delta]/\big(\frac{\omega_n}{\omega}\Lambda[\Delta] \cap  \rho\Lambda[\Delta]\big)
=\rho\Lambda[\Delta]/\rho\frac{\omega_n}{\omega}\Lambda[\Delta]$,}\smallskip

\noindent puisque les $\varphi$-composantes $\rho_\varphi$ de $\rho$ sont étrangères aux polynômes cyclotomiques $\omega_n$. Ainsi:

\begin{Prop}\label{UNormes}
Toujours sous les hypothèses du Théorème \ref{Ucirc}, les sous-groupes normiques attachés aux groupes d'unités et d'unités circulaires sont donnés par les isomorphismes:\smallskip
\begin{itemize}
\item[(i)] $\E_{K_n}^\nu=\wE_{K_n}^{\nu\,\omega}=\varepsilon_n^{\Lambda[\Delta]}\simeq \Lambda[\Delta]/\frac{\omega_n}{\omega}\Lambda[\Delta]$,\smallskip
\item[(ii)] $\E_{K_n}^{\circ\nu}=\wE_{K_n}^{\nu\,\rho\omega}=\varepsilon_n^{\rho\Lambda[\Delta]}
\simeq \rho\Lambda[\Delta]/\rho\frac{\omega_n}{\omega}\Lambda[\Delta]$,\smallskip
\end{itemize}
où  $\varepsilon_n=N_n(\eta^{\,\omega})$ est construit à partir d'une $\Lambda[\Delta]$-base $\eta=(\eta_n)_n$ de $ \,\overset{_\leftarrow}{\E}{}'_K=\varprojlim \,\E'_{K_n}$.
\end{Prop}

\newpage

%%%%%%%%%%%%%%%%%%%%%%%%%%%%%%%%%%%%%%%%%%%%%%%%%%%%%%%%%%%
\subsection{Conjecture faible dans le cas abélien semi-simple $\ell$-décomposé}\medskip
%%%%%%%%%%%%%%%%%%%%%%%%%%%%%%%%%%%%%%%%%%%%%%%%%%%%%%%%%%%

Il est d'usage d'appeler {\em conjecture faible} l'implication: $\,\F_K=1 \Rightarrow \,\C_K=1$,
(cf. e.g. \cite{Ng1,Ng2}); ce qui revient à étudier la conjecture de Greenberg sous l'hypothèse additionnelle $\,\F_K=1$. 
D'après le Théorème \ref{CrSauvage}, dans le cas abélien semi-simple $\ell$-décomposé cela équivaut à supposer $\,\C^{[\ell]}_K=1$.

\begin{Th}\label{Thé}
Soient $K$ un corps abélien réel de degré $d$ et $\ell\nmid d$ un nombre premier impair complètement décomposé dans $K$. Supposons vérifiée la condition de trivialité sauvage $\,\C^{[\ell]}_K=1$; autrement dit que les sous-groupes $\,\Cl_{K_n}^{[\ell]}$ des $\ell$-groupes de classes des corps $K_n$ respectivement engendrés par les images des premiers au-dessus de $\ell$ sont tous triviaux. Alors:\smallskip
\begin{itemize}
\item[(i)] Les groupes d'unités logarithmiques $\,\wE_{K_n}$ coïncident avec les sous-groupes normiques $\,\wE_{K_n}^\nu$:\par
\centerline{$\,\wE_{K_n}^{\phantom{*}}=\,\wE_{K_n}^\nu=\bigcap_{m\ge n} N_{K_{m\!}/\!K_n}(\E'_{K_m})$.}\smallskip

\item[(ii)] Les $\ell$-adifiés des groupes de $\ell$-unités, d'unités et d'unités circulaires des corps $K_n$ admettent les décompositions directes respectives:\smallskip

\centerline{$
\E'_{K_n} =\, \E_K^{\phantom{*}} \oplus \,\wE_{K_n}^{\phantom{*}}; \qquad
\E^{\phantom{*}}_{K_n} =\, \E_K^{\phantom{*}} \oplus \,\wE_{K_n}^{\,\omega}; \qquad
\E^{\circ}_{K_n} =\, \E_K^{\circ} \oplus \,\wE_{K_n}^{\rho\omega}$.}\smallskip

\item[(iii)] Et les $\ell$-groupes de classes circulaires des $K_n$ s'écrivent comme sommes directes:\smallskip

\centerline{$\Cl_{K_n}^{\,\circ}=\,\Cl_K^{\,\circ} \oplus  \,\wE_{K_n}^{\,\omega}/\,\wE_{K_n}^{\rho\omega} \simeq \,\Cl_K^{\,\circ} \oplus \Lambda[\Delta]/\big(\rho\Lambda[\Delta]+\frac{\omega_n}{\omega}\Lambda[\Delta]\big)$.}
\end{itemize}
\end{Th}

\Preuve Supposons $\,\C^{[\ell]}_K=1$, i.e.  $\,\Cl^{[\ell]}_{K_n}=1$ pour tout $n\in\NN$, puisque les places sauvages sont totalement ramifiées dans la tour. Rappelons que cela entraîne $\,\F^{\phantom{*}}_K=\,\F'_K=1$, d'après le Théorème \ref{CrSauvage}; donc, en particulier, l'égalité $\,\wE_K^{\phantom{*}} = \wE_K^\nu$, en vertu de la Proposition \ref{Kuz}; d'où $(i)$.\smallskip

Pour établir $(ii)$, partons de la Proposition \ref{Ram} qui nous donne l'égalité: $\, \E'_K=\,\wE_K^\nu\,\E^{\phantom{*}}_K$ et observons que cette décomposition est directe, comme établi dans la preuve du Théorème \ref{CrSauvage}, i.e.\smallskip

\centerline{ $\, \E'_K=  \,\E^{\phantom{*}}_K \oplus \,\wE_K^\nu =  \,\E^{\phantom{*}}_K \oplus \,\wE_K^{\phantom{*}}$.}\smallskip

Cette identité nous dit en particulier que l'image $\Pl^{[\ell]}_K$ de  $\, \E'_K$ dans le groupe $\wDl^{[\ell]}_K$ des diviseurs logarithmiques de degré nul construits sur les places au-dessus de $\ell$ provient du sous-module des unités $\,\wE_K^{\phantom{*}}$. 
Or, les groupes $\wDl^{[\ell]}_{K_n}$ sont constants, puisque les places au-dessus de $\ell$ sont logarithmiquement inertes dans la $\Zl$-tour cyclotomique $K_\infty/K$; et il en va donc de même des sous-groupes principaux $\Pl^{[\ell]}_{K_n}$, puisque par hypothèse il n'y a pas de capitulation. L'image $\Pl^{[\ell]}_{K_n}$ de $\E'_{K_n}$ dans  $\wDl^{[\ell]}_{K_n}$ provient donc, elle aussi, de $\,\wE_K^{\phantom{*}}$. 
Autrement dit, $\E'_{K_n}$ est engendré par  $\,\E_K^{\phantom{*}}$ et  $\,\wE_{K_n}^{\phantom{*}}$. Et cette décomposition est directe, puisqu'il vient immédiatement, comme plus haut:\smallskip

\centerline{$\E_K^{\phantom{*}} \cap \,\wE_{K_n}^{\phantom{*}} = \,\E_K^{\phantom{*}} \cap \,\wE_{K_n}^{\,\Gamma} = \,\E_K^{\phantom{*}} \cap \,\wE_K^{\phantom{*}}=1$; \quad et donc: \quad $\,\E'_{K_n} =\, \E_K^{\phantom{*}} \oplus \,\wE_{K_n}^{\phantom{*}}$.}\smallskip

Par restriction aux unités, observant que les unités qui sont normes de $\ell$-unités sont en fait normes d'unités, puisque les places sauvages sont totalement ramifiées dans la tour, nous avons:\smallskip

\centerline{$\E_{K_n}^{\phantom{*}} = \E_K^{\phantom{*}} \oplus \big(\,\wE_{K_n}^{\phantom{*}} \cap\,\E_{K_n}^{\phantom{*}}\big) = \,\E_K^{\phantom{*}} \oplus  \,\E_{K_n}^{\nu} = \,\E_K^{\phantom{*}} \oplus \,\wE_{K_n}^{\,\omega}$,}\smallskip

\noindent compte tenu de la Proposition \ref{UNormes}, qui nous donne l'égalité:  $\,\E_{K_n}^\nu=\,\wE_{K_n}^{\nu\,\omega}=\,\wE_{K_n}^{\,\omega}$.\smallskip

Enfin, le cas des unités circulaires se traite de même au moyen de l'isomorphisme canonique\smallskip

\centerline{$\E^{\,\circ}_{K_n}/\E^{\,\circ\,\nu}_{K_n} \simeq \,\E^{\,\circ}_{K}/\E^{\,\circ\,\nu}_K$}\smallskip

\noindent (cf. \cite{Ng1}, Lem. 4.5 ou \cite{Ng2}, Lem. 1.1) et de la Proposition \ref{UNormes}, qui nous donne:  $\,\E_{K_n}^{\,\circ\,\nu}=\,\wE_{K_n}^{\rho\omega}$.\smallskip

L'assertion $(iii)$ s'obtient alors par passage au quotient à partir des deux dernières identités de $(ii)$ et de la descrition donnée par la Proposition \ref{UNormes}.\medskip

\Remarque Les décompositions $\,\E^{\phantom{*}}_{K_n} =\, \E_K^{\phantom{*}} \oplus \,\E_{K_n}^{\nu}$ et 
$\,\E^{\circ}_{K_n} =\, \E_K^{\circ} \oplus \,\E_{K_n}^{\circ\nu}$ données par le Théorème \ref{Thé} font apparaître les groupes d'unités et les sous-groupes circulaires comme sommes respectives des sous-groupes invariants par  $G_n=\Gal(K_n/K)$ et des noyaux de la norme $N_{K_n/K}$. Il suit:\smallskip

%L'inclusion immédiate $N_{K_n/K}(\,\E^\nu_{K_n}) = \,\E_K^\nu \subset \,\wE_K \cap \,\E_F =1$ montre que, sous les hypothèses du Théorème \ref{Thé}, le groupe $\,\E^\nu_{K_n}$ et son sous-groupe circulaire $\,\E^{\,\circ\,\nu}_{K_n}$ sont contenus dans le noyau de la norme $N_{K_n/K}$. Notant $G_n=\Gal(K_n/K)$, nous avons donc:\smallskip

\centerline{$H^1(G_n,\,\E^\nu_{K_n}) = \,\E^\nu_{K_n} /\,\E^{\nu\omega}_{K_n}$;\quad et semblablement: \quad $H^1(G_n,\,\E^{\circ\nu}_{K_n}) = \,\E^{\circ\nu}_{K_n} /\,\E^{\circ\nu\omega}_{K_n}$.}\smallskip

Il résulte alors de la description explicite des numérateurs donnée par la Proposition \ref{UNormes} que le morphisme naturel de $H^1(G_n,\,\E^{\circ\nu}_{K_n}) $ vers $H^1(G_n,\,\E^\nu_{K_n})$ n'est bijectif que si l'élément $\rho$ vaut 1, ce qui équivaut précisément à la conjecture de Greenberg.
En particulier, le Lemme 2.8 de \cite{Ng2} qui affirme sans preuve cette bijectivité revient à postuler la conjecture de Greenberg.

\newpage
%%%%%%%%%%%%%%%%%%%%%%%%%%%%%%%%%%%%%%%%%%%%%%%%%%%%%%%%%%%
\section*{\sc Index des principales notations}
%%%%%%%%%%%%%%%%%%%%%%%%%%%%%%%%%%%%%%%%%%%%%%%%%%%%%%%%%%%
\addcontentsline{toc}{section}{Index des principales notations}

\noindent{\bf Corps et extensions}%%%%%%%%%%%%%%%%%%%%%%%%%%%%%%%%%%%%%%

$K$  : un corps de nombres (totalement réel); et $K_\p$: le complété de $K$ en la place $\p$;\
	
$K_{\infty} = \cup_{n \in \mathbb N} K_n$ avec $[K_n :K]= \ell^n$: la $\Zl$-extension cyclotomique de $K$;\

$K_n$  : le $n$-ième étage de la tour: et $K_{\p_n}$: le complété de $K_n$ en la place $\p_n$;\

$K_n^{nr}$ : la $\ell$-extension abélienne non-ramifiée $\infty$-décomposée maximale de $K_n$;\

$K_n^{\ell d}$ : la $\ell$-extension abélienne non-ramifiée $\ell\infty$-décomposée maximale de $K_n$;\

$K_n^{lc}$ : la pro-$\ell$-extension abélienne localement cyclotomique maximale de $K_n$;\

$K_n^{bp}$ : la pro-$\ell$-extension de Bertrandias-Payan associée à $K_n$;\

$K_n^{\ell r}$ : la pro-$\ell$-extension abélienne $\ell$-ramifiée $\infty$-décomposée maximale de $K_n$;\

$K_\infty^{cd}$ : la pro-$\ell$-extension abélienne partout complètement décomposée maximale de $K_\infty$;\

$K_\infty^{nr}$ : la pro-$\ell$-extension abélienne non-ramifiée $\infty$-décomposée maximale de $K_n$;\

$K_\infty^{\ell r}$ : la pro-$\ell$-extension abélienne $\ell$-ramifiée $\infty$-décomposée maximale de $K_\infty$.

\smallskip

\noindent{\bf Groupes d'idèles} %%%%%%%%%%%%%%%%%%%%%%%%%%%%%%%%%%%%%%

$\R_{K_\p}^{\phantom{*}} = \varprojlim K^\times_\p/K^{\times \ell^m}_\p\!\!$: le compactifié $\ell$-adique du groupe multiplicatif $K^\times_\p$;\

$\U_{K_\p}^{\phantom{*}}$: le sous-groupe unité et $\wU_{K_\p}^{\phantom{*}}$ le groupe des unités logarithmiques dans $\R_{K_\p}^{\phantom{*}}$;\

$\R_{K_\ell}^{\phantom{*}}=\prod_{\l\mid\ell} \R_{K_\l}^{\phantom{*}}$; puis $\,\U_{K_\ell}^{\phantom{*}}=\prod_{\l\mid\ell} \U_{K_\l}^{\phantom{*}}$; et $\,\wU_{K_\ell}^{\phantom{*}}=\prod_{\l\mid\ell} \wU_{K_\l}^{\phantom{*}}$;\

$\mu_{K_\p}$: le $\ell$-groupe des racines de l'unité dans $\R_{K_\p}^{\phantom{*}}$; puis $\mu_{K_\ell}^{\phantom{*}}=\prod_{\l\mid\ell} \mu_{K_\l}^{\phantom{*}}$; et $\mu_{K_\infty}^{\phantom{*}}=\prod_{\p\mid\infty} \mu_{K_\p}^{\phantom{*}}$;\

$\J_K^{\phantom{*}} = \prod^\mathrm{res}_\p \R_{K_\p}^{\phantom{*}}$ : le $\ell$-adifié du groupe des idèles de $K$;\

$\J_K^{[\ell]} = \R_{K_\ell}^{\phantom{*}}\mu_{K_\infty}^{\phantom{*}}\prod_\p \U_{K_\p}^{\phantom{*}}$ : le sous-groupe des $\ell\infty$-idèles;\

$\U_K^{\phantom{*}}=(\prod_\p \U_{K_\p}^{\phantom{*}})\mu_{K_\infty}^{\phantom{*}}$: le sous-groupe unité (au sens ordinaire) de $\J_K^{\phantom{*}}$;

$\wU_K^{\phantom{*}}=(\prod_\p \wU_{K_\p}^{\phantom{*}})\mu_{K_\infty}^{\phantom{*}}$: le groupe des unités logarithmiques (au sens ordinaire) de $\J_K^{\phantom{*}}$.

\smallskip

\noindent{\bf Groupes globaux} %%%%%%%%%%%%%%%%%%%%%%%%%%%%%%%%%%%%%%

$\R_K^{\phantom{*}}=\Zl\otimes_\ZZ K^\times$: le $\ell$-adifié du groupe multiplicatif du corps $K$;\

$\E'_K=\Zl\otimes_\ZZ E'_K$: le $\ell$-adifié du groupe des $\ell$-unités (au sens ordinaire) de $K$;\

$\E_K^{\phantom{*}}=\Zl\otimes_\ZZ E_K=\R_K\cap\,\U_K$: le $\ell$-adifié du groupe des unités (ordinaires) de $K$;\

$\wE_K^{\phantom{*}}=\R_K^{\phantom{*}}\cap\,\wU_K^{\phantom{*}}$: le $\ell$-groupe des unités logarithmiques de $K$;\

$\E_K^\circ$: le $\ell$-groupe des unités circulaires (de Sinnott) attaché au corps (abélien) $K$;\

$\wE_K^{\,\nu}=\bigcap_{n\in\NN}N_{K_n/K}(\E'_{K_n})=\bigcap_{n\in\NN}N_{K_n/K}(\wE_{K_n}^{\phantom{*}})$: le groupe des normes logarithmiques;\

$\E_K^{\,\nu}=\bigcap_{n\in\NN}N_{K_n/K}(\E^{\phantom{*}}_{K_n})$: le sous-groupe des normes d'unités dans $K_{\infty\!}/K$;\

$\E_K^{\circ\nu}=\bigcap_{n\in\NN}N_{K_n/K}(\E^\circ_{K_n})$: le sous-groupe des normes d'unités circulaires dans $K_{\infty\!}/K$.

\smallskip

\noindent{\bf Groupes de classes}%%%%%%%%%%%%%%%%%%%%%%%%%%%%%%%%%%%%%

$\Cl_K^{\phantom{*}}\simeq \J_K^{\phantom{*}}/\,\U_K^{\phantom{*}}\R_K^{\phantom{*}}$: le $\ell$-groupe des classes d'idéaux (au sens ordinaire) du corps $K$;\

$\C ap^{\phantom{*}}_K$: le sous-groupe de capitulation, noyau du morphisme d'extension $\,\Cl_K\rightarrow\,\Cl_{K_\infty}$;\

$\wCl_K^{\phantom{*}}=\wJ_K^{\phantom{*}}/\wU_K^{\phantom{*}}\R_K^{\phantom{*}}$: le $\ell$-groupe des classes logarithmiques  (au sens ordinaire) du corps $K$;\

$\Cl'_K\simeq \J_K/\J_K^{[\ell]} \R_K$: le $\ell$-groupe des $\ell$-classes d'idéaux du corps $K$;\

$\C ap'_K$: le sous-groupe de capitulation, noyau du morphisme d'extension $\,\Cl'_K\rightarrow\,\Cl'_{K_\infty}$;\

$\Cl^{[\ell]}_K$: le sous-groupe de $\,\Cl_K$ engendré par les classes des idéaux premiers au-dessus de $\ell$;\

$\Cl^{\,\circ}_K=\E_K/\E^\circ_K$: le groupe des classes circulaires du corps $K$;\

$\C^\circ_K=\varprojlim\, \Cl^{\,\circ}_{K_n}$: la limite projective des groupes de classes circulaires dans $K_{\infty\!}/K$.

\smallskip

\noindent{\bf Groupes de Galois}%%%%%%%%%%%%%%%%%%%%%%%%%%%%%%%%%%%%%

$\Gamma \ $: le groupe $\Gal(K_{\infty\!}/K)$; $\gamma$ un générateur topologique et $\Lambda=\Zl[[\gamma-1]]$ l'algèbre d'Iwasawa;

$\Gamma_{\!n}=\gamma^{\ell^n\Zl}\ $: le sous-groupe $\Gal(K_{\infty\!}/K_n)$; et $\omega_n=\gamma^{\ell^n}-1$; nous notons $\omega=\gamma-1=\omega_0$;

$\C_K\,=\varprojlim\, \Cl_{K_n} \simeq \Gal(K_\infty^{nr}/K_\infty)$; puis $\C^{[\ell]}_K=\varprojlim \Cl_{K_n}^{[\ell]}$ le sous-module sauvage;\

$\F_K=\varprojlim\, \C ap_{K_n}$: le plus grand sous-module fini du $\Lambda$-module $\,\C_K$;\

$\C'_K\,=\varprojlim\, \Cl'_{K_n} \simeq \Gal(K_\infty^{cd}/K_\infty)$: le module de Kuz'min-Tate attaché au corps $K$;\

$\F'_K=\varprojlim\, \C ap'_{K_n}$: le plus grand sous-module fini du $\Lambda$-module $\,\C'_K$;\

$\T^{bp}_K=\Gal(K^{bp}/K_\infty)$ le groupe de torsion de Bertrandias-Payan.

\newpage
%%%%%%%%%%%%%%%%%%%%%%%%%%%%%%%%%%%%%%%%
%REFERENCES
%%%%%%%%%%%%%%%%%%%%%%%%%%%%%%%%%%%%%%%%

\def\refname{\normalsize{\sc  Références}}

\medskip\noindent
{\small
\begin{tabular}{l}
%{Jean-François {\sc Jaulent}}\\
Institut de Mathématiques de Bordeaux,
Université de {\sc Bordeaux} \& CNRS \\
351 cours de la libération,
F-33405 {\sc Talence} Cedex\\
courriel : Jean-Francois.Jaulent@math.u-bordeaux.fr \quad
\url{https://www.math.u-bordeaux.fr/~jjaulent/}
\end{tabular}
}

 \end{document}